\documentclass[10pt]{article}

\usepackage{amscd,amsmath, amssymb, fancyhdr, graphics, url}

\numberwithin{equation}{section}

\newcommand{\version}{version 9.2,\ \ March 22, 2013}

\makeatletter
\def\x@arrow{\DOTSB\Relbar}
\def\xlongrightarrowfill@{\arrowfill@\relbar\relbar\longrightarrow}
\newcommand{\xlongrightarrow}[2][]{%
        \ext@arrow 0099\xlongrightarrowfill@{#1}{#2}}
\makeatother

\def\eqref#1{(\ref{#1})}
\newcommand{\goth}{\mathfrak}

\newcommand{\arrow}{{\:\longrightarrow\:}}
\newcommand{\Z}{{\Bbb Z}}
\def\C{{\Bbb C}}
\newcommand{\R}{{\Bbb R}}
\newcommand{\Q}{{\Bbb Q}}

\newcommand{\6}{\partial}
\def\1{\sqrt{-1}\:}
\newcommand{\restrict}[1]{{\left|_{{\phantom{|}\!\!}_{#1}}\right.}}
\newcommand{\cntrct}                
{\hspace{2pt}\raisebox{1pt}{\text{$\lrcorner$}}\hspace{2pt}}

\renewcommand{\tilde}{\widetilde}
\renewcommand{\bar}{\overline}
\renewcommand{\phi}{\varphi}
\renewcommand{\epsilon}{\varepsilon}
\renewcommand{\geq}{\geqslant}
\renewcommand{\leq}{\leqslant}

\newcommand{\Teich}{\operatorname{\sf Teich}}
\newcommand{\Comp}{\operatorname{\sf Comp}}
\newcommand{\NS}{\operatorname{\sf NS}}
\newcommand{\Per}{\operatorname{\sf Per}}
\newcommand{\Perspace}{\operatorname{{\Bbb P}\sf er}}

\newcommand{\Id}{\operatorname{Id}}

\newcommand{\Hom}{\operatorname{Hom}}
\newcommand{\Sym}{\operatorname{Sym}}
\newcommand{\Pic}{\operatorname{Pic}}

\newcommand{\Aut}{\operatorname{Aut}}
\newcommand{\Mon}{\operatorname{\sf Mon}}

\newcommand{\Alb}{\operatorname{Alb}}

\newcommand{\Diff}{\operatorname{\sf Diff}}
\newcommand{\Spin}{\operatorname{Spin}}
\newcommand{\codim}{\operatorname{codim}}

\newcommand{\rk}{\operatorname{rk}}

\newcommand{\Spec}{\operatorname{Spec}}

\renewcommand{\Re}{\operatorname{Re}}
\renewcommand{\Im}{\operatorname{Im}}

\newcounter{Mycounter}[section]
\newcounter{lemma}[section]
\renewcommand{\thelemma}{{Lemma \thesection.\arabic{lemma}}}
\newcommand{\lemma}{%
    \setcounter{lemma}{\value{Mycounter}}
    \refstepcounter{lemma}
    \stepcounter{Mycounter}
    {\noindent \bf \thelemma:\ }}

\newcounter{claim}[section]
\renewcommand{\theclaim}{{Claim \thesection.\arabic{claim}}}
\newcommand{\claim}{%
    \setcounter{claim}{\value{Mycounter}}
    \refstepcounter{claim}
    \stepcounter{Mycounter}
    {\noindent \bf \theclaim:\ }}

\newcounter{sublemma}[section]
\newcounter{corollary}[section]
\renewcommand{\thecorollary}{{Corollary \thesection.\arabic{corollary}}}
\newcommand{\corollary}{%
    \setcounter{corollary}{\value{Mycounter}}
    \refstepcounter{corollary}
    \stepcounter{Mycounter}
    {\noindent \bf \thecorollary:\ }}

\newcounter{theorem}[section]
\renewcommand{\thetheorem}{{Theorem \thesection.\arabic{theorem}}}
\newcommand{\theorem}{%
    \setcounter{theorem}{\value{Mycounter}}
    \refstepcounter{theorem}
    \stepcounter{Mycounter}
    {\noindent \bf \thetheorem:\ }}

\newcounter{conjecture}[section]
\newcounter{proposition}[section]
\renewcommand{\theproposition}
      {{Proposition \thesection.\arabic{proposition}}}
\newcommand{\proposition}{%
    \setcounter{proposition}{\value{Mycounter}}
    \refstepcounter{proposition}
    \stepcounter{Mycounter}
    {\noindent \bf \theproposition:\ }}

\newcounter{definition}[section]
\renewcommand{\thedefinition}
      {{Definition~\thesection.\arabic{definition}}}
\newcommand{\definition}{%
    \setcounter{definition}{\value{Mycounter}}
    \refstepcounter{definition}
    \stepcounter{Mycounter}
    {\noindent \bf \thedefinition:\ }}

\newcounter{example}[section]
\renewcommand{\theexample}{{Example \thesection.\arabic{example}}}
\newcommand{\example}{%
    \setcounter{example}{\value{Mycounter}}
    \refstepcounter{example}
    \stepcounter{Mycounter}
    {\noindent \bf \theexample:\ }}

\newcounter{remark}[section]
\renewcommand{\theremark}{{Remark \thesection.\arabic{remark}}}
\newcommand{\remark}{%
    \setcounter{remark}{\value{Mycounter}}
    \refstepcounter{remark}
    \stepcounter{Mycounter}
    {\noindent \bf \theremark:\ }}

\newcounter{problem}[section]
\newcounter{question}[section]
\makeatletter

\@addtoreset{equation}{section} \@addtoreset{footnote}{section}
\makeatother

\def\blacksquare{\hbox{\vrule width 5pt height 5pt depth 0pt}}
\def\endproof{\blacksquare}

\begin{document}

\begin{center}
{\LARGE\bf
Mapping class group\\[1mm] and a global Torelli theorem\\[1mm] for hyperk\"ahler manifolds\\[4mm]
}

 Misha
Verbitsky\footnote{Partially supported by RFBR grants
 12-01-00944-Á,  NRI-HSE 
Academic Fund Program in 2013-2014, research grant
12-01-0179, and
AG Laboratory NRI-HSE, RF government grant, ag. 11.G34.31.0023.


}

\end{center}

{\small \hspace{0.15\linewidth}
\begin{minipage}[t]{0.7\linewidth}
{\bf Abstract} \\
{\bf A mapping class group} of an oriented manifold
is a quotient of its diffeomorphism group by the  isotopies.
We compute a mapping class group of a hypek\"ahler
manifold $M$, showing that it is commensurable to
an arithmetic lattice in $SO(3, b_2-3)$. A Teichm\"uller
space of $M$ is a space of complex structures on $M$ up
to isotopies. We define  {\bf a birational Teichm\"uller
space} by identifying certain points corresponding
to bimeromorphically equivalent manifolds. We show that
the period map gives the isomorphism between connected
components of the birational Teichm\"uller space and the
corresponding period space 
$SO(b_2-3, 3)/SO(2)\times SO(b_2 -3, 1)$. We use this
result to obtain a Torelli theorem identifying each
connected component of  the birational moduli space with 
a quotient of a period space by an arithmetic group.
When $M$ is a Hilbert scheme 
of $n$ points on a K3 surface, with $n-1$
a prime power, our Torelli theorem implies the usual 
Hodge-theoretic birational Torelli theorem
(for other examples of hyperk\"ahler manifolds,
the Hodge-theoretic Torelli theorem is known
to be false).
\end{minipage}
}

\tableofcontents


\section{Introduction}


\subsection{Hyperk\"ahler manifolds and their moduli}

Throughout this paper, {\bf a hyperk\"ahler manifold}
is a compact, holomorphically symplectic manifold
of K\"ahler type, simply connected and with
$H^{2,0}(M)=\C$. In the literature, such
manifolds are often called {\bf simple}, or {\bf irreducible}.
For an explanation of this term and an
introduction to hyperk\"ahler structures, please 
see Subsection \ref{_hk_stru_Subsection_}.

We shall say that a complex structure $I$ on $M$
is {\bf of hyperk\"ahler type} if $(M,I)$ is
a hyperk\"ahler manifold. 

There are many different ways to define the moduli
of complex structures. In this paper we use the
earliest one, which is due to Kodaira-Spencer and Kuranishi.
Let $M$ be an oriented manifold, ${\goth I}$
the space of all complex structures of 
hyperk\"ahler type, compatible with orientation,
and ${\cal M}:={\goth I}/\Diff$ 
its quotient by the group of oriented 
diffeomorphisms.\footnote{Throughout this paper,
we speak of oriented diffeomorphisms, but the reasons
for this assumption are purely historical. We could
omit the mention of orientation, and most of the results will
remain valid.}
We call ${\cal M}$ {\bf the moduli space} of complex structures 
of hyperk\"ahler type (or just ``the moduli space'') of
$M$. This space is usually non-Hausdorff. 

For a hyperk\"ahler manifold, the non-Hausdorff points of $\cal M$ are 
easy to control, due to a theorem of D. Huybrechts
(\ref{_bira_points_inseparable_Theorem_}).
If $I_1, I_2\in {\cal M}$ are inseparable points
in ${\cal M}$, then the corresponding hyperk\"ahler
manifolds are bimeromorphic 
(\ref{_non_Hausdorff_marked_and_Teich_Proposition_}).

In many cases, the moduli of complex structures 
on $M$ can be described in terms of Hodge structures
on the cohomology of $M$. Such results are called 
{\em Torelli theorems}. In this paper, we state
a Torelli theorem for hyperk\"ahler manifolds,
using the language of mapping class group
and Teichm\"uller spaces.

This approach to the Torelli-type problems
was pioneered by A. Todorov in several
important preprints and papers (\cite{_Todorov:torsion_}, 
\cite{_Todorov:WP_}; see also \cite{_LTYZ:moduli_}).

\subsection{Teichm\"uller space of a hyperk\"ahler manifold}

To define the period space for hyperk\"ahler manifolds,
one uses the so-called Bogomolov-Beauville-Fujiki (BBF) form
on the second cohomology. Historically, it was the
BBF form which was defined in terms of the period space,
and not vice versa, but the other way around is more
convenient.

\hfill

Let $\Omega$ be a holomorphic symplectic form on $M$.
Bogomolov and Beauville
(\cite{_Bogomolov:defo_}, \cite{_Beauville_}) defined the
following bilinear symmetric 2-form on $H^2(M)$:
\begin{equation}\label{_BBF_form_on_H^11_Equation_}
\begin{aligned}
\tilde q(\eta,\eta'):= &
  2\int_M \eta\wedge\eta'  \wedge \Omega^{n-1}
   \wedge \bar \Omega^{n-1} -\\
 &-\frac{n-1}{n}\frac{\left(\int_M \eta \wedge \Omega^{n-1}\wedge \bar
   \Omega^{n}\right) \left(\int_M \eta' \wedge
   \Omega^{n}\wedge \bar \Omega^{n-1}\right)}{\int_M \Omega^{n}
   \wedge \bar \Omega^{n}}
\end{aligned}
\end{equation}
where $n=\dim_{\Bbb H}M$.

\hfill

\remark 
The form $\tilde q$ is compatible with the Hodge decomposition,
which is seen immediately from its definition. Also,
$\tilde q(\Omega, \bar \Omega)>0$.

\hfill

The form $\tilde q$ is topological by its nature.

\hfill

\theorem\label{_BBF_Theorem_}
(\cite{_Fujiki:HK_})
Let be a simple hyperk\"ahler
manifold of real dimension $4n$. Then there exists a
bilinear, symmetric, primitive non-degenerate
integral 2-form $q:\; H^2(M,\Z)\otimes H^2(M,\Z)\arrow \Z$
and a constant $c\in \Z$
such that 
\begin{equation}\label{_BBF_form_Equation_}
\int_M \eta^{2n}= cq(\eta,\eta)^n,
\end{equation}
for all $\eta\in H^2(M)$. Moreover,
$q$ is proportional to the form $\tilde q$ of
\eqref{_BBF_form_on_H^11_Equation_}, and has signature
$(+,+,+, -,-,-, ...)$.

\endproof

\hfill

\remark
If $n$ is odd, the equation \eqref{_BBF_form_Equation_}
determines $q$ uniquely, otherwise -- up to a sign.
To choose a sign, we use \eqref{_BBF_form_on_H^11_Equation_}.

\hfill

\definition\label{_BBF_Definition_}
Let $M$ be a hyperk\"ahler manifold,
and $\Omega$ a holomorphic symplectic form 
on $M$. {\bf Beauville-Bogomolov-Fujiki form} on $M$
is a form $q:\; H^2(M,\Q)\otimes H^2(M,\Q)\arrow \Q$
which satisfies \eqref{_BBF_form_Equation_},
and has $q(\Omega, \bar\Omega)>0$.

\hfill

\definition
Let $(M,I)$ be a compact hyperk\"ahler manifold, ${\goth I}$
the set of oriented complex structures
of hyperk\"ahler type on $M$, and $\Diff_0(M)$ the group
of isotopies. The quotient space
$\Teich:={\goth I}/\Diff_0(M)$ is called 
{\bf the Teichm\"uller space} of $(M, I)$,
and the quotient of $\Teich$ over a whole oriented
diffeomorphism group {\bf the coarse moduli space of $(M,I)$}.

\hfill

\remark
In a similar way one defines the moduli of K\"ahler structures
or of complex structures on a given K\"ahler or complex manifold.
This approach was originaly suggested by Kodaira and Spencer
in their fundamental work on deformation theory 
(\cite{_Kodaira_Spencer:I_}, \cite{_Kodaira_Spencer:III_}).
Kodaira and Spencer constructed a local moduli space for
the special cases when obstructions to deformation vanish.
The results of Kodaira-Spencer were obtained in full generality
by M. Kuranishi, who constructed a finite-dimensional
complex-analytic slice of the action of diffeomorphism
group on the space of all integrable almost
complex structures (\cite{_Kuranishi:new_}, 
\cite{_Kuranishi:note_}).
The precise correspondence between 
the results of Kuranishi and the action of
the diffeomorphism group was spelled out by
A. Douady in his Bourbaki talk, \cite{_Douady:Bourbaki_}.

\hfill

\remark
As shown by F. Catanese \cite[Proposition 15]{_Catanese:moduli_},
for K\"ahler manifolds with trivial canonical bundle,
e.g. for the hyperk\"ahler manifolds,
the Teichm\"uller space is locally isomorphic
to the Kuranishi moduli space.

\hfill

\definition
Let $(M,I)$ be a simple hyperk\"ahler manifold,
and $\Teich$ its Teichm\"uller space. For any $J\in \Teich$,
$(M,J)$ is also a simple hyperk\"ahler manifold, as seen
from \ref{_irredu_topolo_Lemma_} below, hence
$H^{2,0}(M,J)$ is one-dimensional. Consider a
map $\Per:\; \Teich \arrow {\Bbb P}H^2(M, \C)$,
sending $J$ to the line $H^{2,0}(M,J)\in {\Bbb P}H^2(M, \C)$.
Clearly, $\Per$ maps $\Teich$ into an open subset of the 
quadric, defined by
\begin{equation}\label{_perspace_lines_Equation_}
\Perspace:= \{l\in {\Bbb P}H^2(M, \C)\ \ | \ \  q(l,l)=0, q(l, \bar l) >0\}.
\end{equation}
The map $\Per:\; \Teich \arrow \Perspace$ is called {\bf the period map},
and the set $\Perspace$ {\bf the period space}.

\hfill

The following fundamental theorem is due
to F. Bogomolov \cite{_Bogomolov:defo_}.

\hfill

\theorem\label{_bogomo_smooth_Theorem_}
(\cite{_Bogomolov:defo_}) Let $M$ be a simple hyperk\"ahler manifold,
and $\Teich$ its Teichm\"uller space. Then
the period map $\Per:\; \Teich \arrow \Perspace$ is 
locally an unramified covering (that is, an etale map).
\endproof

\hfill

\remark
Bogomolov's theorem implies that $\Teich$ is smooth.
However, it is not necessarily Hausdorff
(and it is non-Hausdorff even in the 
simplest examples).

\hfill

\remark
D. Huybrechts has shown that $\Per$ is surjective
(\cite{_Huybrechts:basic_},  Theorem 8.1).

\hfill

\remark\label{_teich_finite_compo_Remark_}
Using the boundedness results of Kollar and
Matsusaka (\cite{_Kollar_Matsusaka:finiteness_}), 
D. Huybrechts has shown that 
the space $\Teich$ has only a 
finite number of connected components
(\cite{_Huybrechts:finiteness_}, Theorem 2.1).

\hfill

The moduli space ${\cal M}$ of complex structures 
of hyperk\"ahler type on $M$ is 
a quotient of $\Teich$ by the action
of the mapping class group $\Gamma:=\Diff/\Diff_0$
of diffeomorphisms up to isotopies. There is an 
interesting intermediate group $\Diff_H$ of all
diffeomorphisms acting trivially on $H^2(M)$.
One has $\Diff_0 \subset \Diff_H \subset \Diff$.
The corresponding quotient $\Teich/\Diff_H$ is called
{\bf the coarse, marked moduli space} of complex structures,
and its points -- {\bf marked hyperk\"ahler manifolds}. 
To choose a marking it means to choose a basis
in the cohomology of $M$. The period map
is well defined on $\Teich/\Diff_H$. 

We don't use the marked moduli space in this paper,
because the Teichm\"uller space serves the same purpose.
In the literature on moduli spaces, the marked moduli space 
is used throughout, but these results are easy to 
translate to the Teichm\"uller spaces'
language using the known facts about the
mapping class group.

For a K3 surface, the Teichm\"uller space is not
Hausdorff. However, a quotient of polarized moduli space
by the mapping class group {\em is}
Hausdorff  and quasi-projective (\cite{_Viehweg:moduli_}).
Moreover, a version of Torelli theorem is valid,
providing an isomorphism between $\Teich/\Gamma$
and $\Perspace/O^+(H^2(M, \Z))$.\footnote{For an explanation
of $O^+$, please see \ref{_spinorial_norm_Definition_}.}
This result has a long
history, with many people contributing to different
sides of the picture, but its conclusion could be found in 
\cite{_Burns_Rappoport_} and \cite{_Siu:K3_}.

One could state this Torelli theorem as 
a result about the Hodge structures,
as follows. The Torelli theorem 
claims that there  is a bijective correspondence
between isomorphism classes of K3 surfaces and the
set of isomorphism classes of appropriate Hodge structures 
on a 22-dimensional space equipped with an integer lattice,
a spin orientation (\ref{_spin_orienta_hk_Remark_}) 
and an integer quadratic form.

It is natural to expect that this last result would be
generalized to other hyperk\"ahler manifolds, but such a
straightforward generalization is invalid. 
In \cite{_Debarre:Torelli_}, O. Debarre has shown that
there exist birational hyperk\"ahler manifolds which are 
non-isomorphic, but have the same periods.
A hope to have a Hodge theoretic Torelli theorem
for birational moduli was extinguished in early
2000-ies. As shown by Yo. Namikawa in a beautiful
(and very short) paper \cite{_Namikawa:Torelli_},
there exist hyperk\"ahler manifolds $M, M'$
which are not bimeromorphically equivalent,
but their second cohomology have equivalent
Hodge structures.

For the benefit of the reader, we give here a brief
reprise of the
Na\-mi\-ka\-wa's construction.  Let $T$ be a compact,
complex, 2-dimensional torus, and $T^{[n]}$
its Hilbert scheme. The Albanese map
$T^{[n]}\stackrel \Alb\arrow T$ is a locally trivial
fibration. Denote by $T^{[n]}_T$
{\bf the generalized Kummer manifold}, 
 $T^{[n]}_T:= \Alb^{-1}(0)$. When
$n=2$, it is a K3 surface obtained
from the torus using the Kummer construction. 
For $n>2$, the Hodge structure 
on $H^2(T^{[n]}_T)$ is easy to describe.
One has
\[
H^2(T^{[n]}_T)\cong \Sym^2(H^1(T))\oplus \R \eta,
\]
where $\eta$ is the fundamental class of the 
exceptional divisor of $M:=T^{[n]}_T$.
Therefore, $H^2(M)$ has the same Hodge structure as 
$M' =(T^*)^{[n]}_{T^*}$, where $T^*$ is the dual torus.
However, the manifolds $M$ and $M'$ are not
bimeromorphically equivalent, when $T$ is generic. 
This is easy to see, for instance, for $n=3$, because the 
exceptional divisor of $M=T^{[3]}_T$ is a
trivial $\C P^1$-fibration over $T$,
and the exceptional divisor of $M'=(T^*)^{[3]}_{T^*}$
is fibered over $T^*$ likewise. Since 
bimeromorphic maps of holomorphic symplectic
varieties are non-singular in codimension 2,
any bimeromorphic isomorphism between $M$ and $M'$
would bring a bimeromorphic isomorphism between
these divisors, and therefore between $T$ and $T^*$,
which is impossible for general $T$.

A less elementary construction,
due to E. Markman, gives a counterexample 
to the Hodge-theoretic global Torelli theorem 
when $M=K3^{[n]}$ is the Hilbert scheme
of points on a K3 surface, and $n-1$
is not a prime power (\cite{_Markman:constra_}). When $n-1$ 
is a prime power, a Hodge-theoretic birational
Torelli theorem holds true (Subsection 
\ref{_Hodge_Torelli_Subsection_}).

We are going to prove a different version of Torelli theorem,
using the language of Teichm\"uller spaces and the mapping
class groups.

\subsection{The birational Teichm\"uller space}

The Teichm\"uller space approach allows one to state the
Torelli theorem for hyperk\"ahler manifolds as it is done
for curves. However, before any theorems can be stated, we
need to resolve the issue of non-Hausdorff points.

\hfill

\definition\label{_insepara_Definition_}
Let $M$ be a topological space. We say that points
$x, y\in M$ are {\bf inseparable} (denoted $x\sim y$)
if for any open subsets $U\ni x, V\ni y$, one 
has $U \cap V\neq \emptyset$.

\hfill

\remark 
As follows from \ref{_non_Hausdorff_marked_and_Teich_Proposition_}
and \ref{_bira_points_inseparable_Theorem_},
inseparable points on a Teichm\"uller space correspond to
bimeromorphically equivalent hy\-per\-k\"ah\-ler manifolds. 

\hfill

\theorem
Let $\Teich$ be a Teichm\"uller space of a hyperk\"ahler
manifold, and $\sim$ the inseparability relation defined above.
Then $\sim$ is an equivalence relation. Moreover,
the quotient $\Teich_b:=\Teich\!/{}_\sim$ is a smooth, Hausdorff
complex analytic manifold.

{\bf Proof:} \ref{_Hausdorff_redu_Theorem_}, 
\ref{_Teich_b_manifold_Theorem_}. \endproof

\hfill

We call the quotient $\Teich\!/{}_\sim$ 
{\bf the birational Teichm\"uller space},
denoting it as $\Teich_b$. The operation
of taking the quotient $.../{}_\sim$ as above
has good properties in many situations, and
brings similar results quite often.
We call $W\!/{}_\sim$
{\bf the Hausdorff reduction} of $W$
whenever it is Hausdorff
(see Subsection \ref{_Hausdorff_red_Subsection_} 
for a detailed expos\'e).

\subsection{The mapping class group of a hyperk\"ahler manifold}

Define the mapping class group 
$\Gamma:= \Diff/\Diff_0$ of a manifold $M$
as a quotient of the group of oriented diffeomorphisms of $M$ 
by isotopies. Clearly, $\Gamma$ acts on $H^2(M, \R)$
perserving the integral structure. We are able to determine
the group $\Gamma$ up to commensurability, proving that
it is commensurable to an arithmetic group
$O(H^2(M, \Z), q)$ of finite covolume in 
$O(3, b_2(M)-3)$. 

\hfill

\theorem \label{_mapping_class_intro_Equation_}
Let $M$ be a compact, simple hyperk\"ahler manifold, and
$\Gamma= \Diff/\Diff_0$ its mapping class group.
Then $\Gamma$ acts on $H^2(M, \R)$ preserving the
Bogomolov-Beauville-Fujiki form. Moreover, the corresponding
homomorphism $\Gamma \arrow O(H^2(M, \Z), q)$ has finite
kernel, and its image has finite index in
$O(H^2(M, \Z), q)$.

{\bf Proof:} This is \ref{_Mapp_cl_group_hk_Theorem_}. 
\endproof

\hfill

Using results of E. Markman (\cite{_Markman:constra_}),
it is possible to compute the mapping class
group for a Hilbert scheme of points on a K3
surface $M=K3^{[n]}$, when $n-1$ is a prime
power (\ref{_mono_ref_Theorem_}).

\subsection{Teichm\"uller space and Torelli-type theorems}

The following version of the Torelli theorem is proven in
Section \ref{_HK_lines_Section_}.

\hfill

\theorem\label{_Torelli_intro_Theorem_}
Let $M$ be a compact,
simple hyperk\"ahler manifold, and $\Teich_b$ its birational
Teichm\"uller space. Consider the period map
$\Per:\; \Teich_b:\; \arrow \Perspace$, where $\Perspace$
is the period space defined as in \eqref{_perspace_lines_Equation_}.
Then $\Per$ is a diffeomorphism, for each connected
component of $\Teich_b$.

{\bf Proof:} This is  \ref{_global_Torelli_Theorem_}. \endproof

\hfill

The proof of \ref{_Torelli_intro_Theorem_} is
obtained by using the quaternionic structures,
associated with holomorphic symplectic structures
by the Calabi-Yau theorem, and the corresponding
rational lines in $\Teich$ and $\Perspace$. 

If one wants to obtain a more traditional
Torelli-type theorem, one should consider the 
set of equivalence classes of complex structures up
to birational equivalence. This set can be
interpreted in terms of the Teichm\"uller 
space as follows. 

Consider the action of the mapping class
group $\Gamma$ on the Teichm\"uller space $\Teich$, and let $\Teich^I$
be a connected component of $\Teich$ containing a given
complex structure $I$. Denote by $\Gamma_I\subset \Gamma$ 
a subgroup of $\Gamma$
preserving $\Teich^I$. Since $\Teich$
has only a finite number of connected components 
(\cite{_Huybrechts:finiteness_}, Theorem 2.1),
$\Gamma_I$ has a finite index in $\Gamma$.
The coarse moduli space of complex structures
on $M$ is $\Teich^I/\Gamma_I$, and the birational
moduli is $\Teich^I_b/\Gamma_I$, where $\Teich_b^I$
is the appropriate connected component of $\Teich_b$. 
\ref{_Torelli_intro_Theorem_} immediately implies
the following Torelli-type result.

\hfill

\theorem\label{_bira_modu__Torelli_Theorem_}
Let $M$ be a compact, simple hyperk\"ahler manifold,
 ${\cal M}_b:= \Teich^I_b/\Gamma_I$ a connected component
of the birational moduli space defined above, and
\begin{equation}\label{_period_bira_modu_Equation_}
{\cal M}_b\stackrel{\Per} \arrow\Perspace/{\Gamma_I}
\end{equation}
the corresponding period map. Then 
\eqref{_period_bira_modu_Equation_}
is a bijection. \endproof

\hfill

\remark
The image $i(\Gamma_I)$ of $\Gamma_I$ in $O(H^2(M, \Z), q)$
has finite index (\ref{_mapping_class_intro_Equation_}).
Therefore, it is an arithmetic sugbroup of finite covolume.

\hfill

Comparing this with \ref{_bira_modu__Torelli_Theorem_},
we immediately obtain the following corollary.

\hfill

\corollary\label{_bira_modu_Corollary_}
Let $M$ be a compact, simple hyperk\"ahler manifold,
and ${\cal M}_b$  a connected 
component of its birational moduli space, obtained as above.
Then ${\cal M}_b$ is isomorphic to a quotient of a homogeneous
space 
\[ 
  \Perspace=\frac{O(b_2-3,3)}{SO(2)\times O(b_2-3, 1)}
\]
by an action of an arithmetic subgroup
$i(\Gamma_I)\subset O(H^2(M, \Z), q)$.\footnote{For 
this interpretation of $\Perspace$, please see 
Subsection \ref{_period_space_grassman_Subsection_}.}
\endproof

\hfill

In a traditional version of Torelli theorem,
one takes a quotient of $\Perspace$
by $O^+(H^2(M, \Z), q)$ instead of 
$i(\Gamma_I)\subset O^+(H^2(M, \Z), q)$.%
\footnote{$O^+(H^2(M, \Z), q)$ is a group of orthogonal
maps with positive spin norms (\ref{_spinorial_norm_Definition_}).}
However, such a result cannot be valid,
as shown by Namikawa. \ref{_bira_modu_Corollary_}
explains why this occurs: for Namikawa's examples, 
the group $i(\Gamma_I)$ is a proper subgroup in
$O^+(H^2(M, \Z), q)$, and the composition
\begin{equation}\label{_moduli_to_Torelli_compo_Equation_}
 {\cal M}_b\arrow \Perspace/\Gamma_I 
 \arrow \Perspace/O^+(H^2(M, \Z), q)
\end{equation}
is a finite quotient map. We obtained the following corollary.

\hfill

\corollary\label{_fini_quo_Corollary_}
Let $M$ be a compact, simple hyperk\"ahler manifold,
 ${\cal M}_b$ a connected component of its birational moduli space,
and 
\begin{equation}\label{_periods_trad_Equation_}
{\cal M}_b\arrow  \Perspace/O^+(H^2(M, \Z), q)
\end{equation}
the corresponding period map. Then
\eqref{_periods_trad_Equation_}
is a finite quotient. \endproof

\hfill

\remark
Please notice that 
the space $\Perspace/O^+(H^2(M, \Z), q)$ is usually 
non-Hausdorff. However, it can be made 
Hausdorff if one introduces additional
structures (such as a polarization), and then
\ref{_fini_quo_Corollary_} becomes more useful.

\hfill

For the Hilbert scheme of $n$ points on a K3 surface,
the image of $\Gamma_I$ in $O^+(H^2(M, \Z))$
was computed by E. Markman in \cite{_Markman:constra_}
(see \ref{_mono_ref_Theorem_}). When $n-1$ is a prime power, 
$i(\Gamma_I)=O^+(H^2(M, \Z))$, and
the composition \eqref{_periods_trad_Equation_}
is an isomorphism, which is used to obtain 
the usual (Hodge-theoretic) version of 
Torelli theorem.

\subsection{A Hodge-theoretic Torelli theorem for
$K3^{[n]}$}
In \cite{_Markman:mono_}, \cite{_Markman:constra_},
E. Markman has proved many vital results on the
way to compute the mapping class group
of a Hilbert scheme of points on K3 (denoted
by $K3^{[n]}$). Markman's starting point was
the notion of a monodromy group of a hyperk\"ahler manifold.
A monodromy group of $M$ is the group generated by 
monodromy of the Gauss-Manin local systems
for all deformations of $M$ (see Subsection
\ref{_mono_hype_Subsection_} 
for a more detailed description).
In Subsection
\ref{_mono_hype_Subsection_}, we relate
the monodromy group $\Mon$ to the mapping class
group $\Gamma_I$, showing that $\Mon$ is isomorphic
to an image of $\Gamma_I$ in $PGL(H^2(M, \C))$.
For $M=K3^{[n]}$,  Markman has computed the
monodromy group, using the action of Fourier-Mukai
transform in the derived category of coherent sheaves. 
He used this computation to show that the
standard (Hodge-theoretic) global Torelli theorem
fails on $K3^{[n]}$, unless $n-1$ is a prime
power. We complete Markman's analysis of
global Torelli problem for $K3^{[n]}$,
proving the following.

\hfill

\theorem
Let $M=K3^{[n]}$ be a Hilbert scheme of points
on a K3 surface, where $n-1$ is a prime power, and
$I_1, I_2$ deformations of complex structures on 
$M$. Assume that the Hodge structures on 
$H^2(M, I_1)$ and $H^2(M, I_2)$ are isomorphic,
and this isomorphism is compatible with the
Bogomolov-Beauville-Fujiki form and the natural
spin orientation on $H^2(M, I_1)$ and $H^2(M, I_2)$.
(\ref{_spin_orienta_hk_Remark_}).
Then $(M, I_1)$ is bimeromorphic to $(M, I_2)$.

\hfill

{\bf Proof:} This is \ref{_Torelli_prime_powe_Theorem_}. \endproof

\hfill

\remark
E. Markman in \cite{_Markman:constra_} constructed
counterexamples to the  Hodge-theoretic
global Torelli problem for $K3^{[n]}$,
where $n-1$ is not a prime power.

\subsection{Moduli of polarized hyperk\"ahler varieties}

For another application of 
\ref{_bira_modu_Corollary_}, fix an integral
class $\eta \in H^2(M, \Z)$, $q(\eta,\eta)>0$,
and let $\Teich_{\pm,\eta}$ be a divisor in the
connected component of the Teichm\"uller space consisting of all
$I$ with $\eta\in H^{1,1}(M,I)$.
For a general $I\in \Teich_{\pm,\eta}$, $\eta$ or $-\eta$ is a K\"ahler
class on $(M,I)$ (\cite{_Huybrechts:cone_}; 
see also \ref{_Kahler_cone_Pic=1_Theorem_}).
However, there could be special points where
$\pm \eta$ is not K\"ahler.

Let $\Teich_\eta^I$ be a connected component of
$\Teich_{\pm,\eta}$, containing $I\in \Teich$,
in such a way that $\eta$ is a K\"ahler class on $(M,I)$. 
Denote by $\bar {\cal M}_\eta$ the quotient of
$\Teich_\eta^I$ by the subgroup $\Gamma_\eta^I$ of the mapping class
group fixing $\eta$ and preserving the component
$\Teich_\eta^I$. The same argument as above
can be used to show that $\Gamma_\eta^I$
is commensurable to an arithmetic subgroup
in $SO(\eta^\bot)$, where $\eta^\bot \subset H^2(M, \R)$
is an orthogonal complement to $\eta$.

We call $\bar {\cal M}_\eta$ 
{\bf a connected component of the moduli space of weakly
polarized hyperk\"ahler manifolds.}
A corresponding component ${\cal M}_\eta$ of the moduli 
of polarized hyperk\"ahler
manifolds  is an open subset of $\bar {\cal M}_\eta$ 
consisting of all $I$ for which $\eta$ is 
K\"ahler. Since $\pm\eta$ is a K\"ahler class whenever
 $\Pic(M)=\langle \eta\rangle$ (\ref{_Kahler_cone_Pic=1_Theorem_}),
${\cal M}_\eta$ is dense in $\bar {\cal M}_\eta$.
It is known (due to the general theory
which goes back to Viehweg, Grothendieck and Kodaira-Spencer)
that ${\cal M}_\eta$ is Hausdorff and quasiprojective
(see e.g. \cite{_Viehweg:moduli_}, \cite{_GHK:moduli_HK_}).

Recall that for any simple hyperk\"ahler manifold $M$,
the space $H^2(M, \R)$ has signature $(+,+,+, -, ..., -)$,
with a fixed orientation on each of the positive
3-planes. 
The period space for weakly polarized hyperk\"ahler
manifolds is defined as
\begin{equation}\label{_weakly_perspace_Equation_}
\Perspace_{\pm,\eta}:= 
\{l\in {\Bbb P}H^2(M, \C)\ \ | \ \  q(l,l)=0, \
q(\eta,l)=0, \ q(l, \bar l) >0\},
\end{equation}
where the 3-plane $\langle \eta, \Re l, \Im l\rangle$
has the same orientation as for a hyperk\"ahler structure.
 Then the corresponding 
period map $\Teich_\eta \arrow \Perspace_\eta$
induces an isomorphism from the 
Hausdorff reduction $\Teich_{\eta,b}^I$ of 
 $\Teich_\eta^I$ to $\Perspace_\eta$,
as follows from \ref{_Torelli_intro_Theorem_}.

Notice that, unless one fixes the orientation of
the 3-plane $\langle \eta, \Re l, \Im l\rangle$,
the space $\{l\in {\Bbb P}H^2(M, \C)\ \ | \ \  q(l,l)=0, \
q(\eta,l)=0, \ q(l, \bar l) >0\}$ would have {\it two} 
connected components.

We define {\bf a connected component of the birational moduli space
of weakly polarized hyperk\"ahler manifolds}
$\bar {\cal M}_{b,\eta}$  as a quotient of the component
$\Teich_{b,\eta}^I$ by the 
corresponding mapping class group $\Gamma_\eta^I$.
It is obtained from $\bar {\cal M}_\eta$
by identifying inseparable points.

Just as in Subsection
\ref{_period_space_grassman_Subsection_},
we may identify the period space
$\Perspace_\eta$ with the Grassmannian of 
positive 2-planes in $\eta^\bot$. This gives
\[
\Perspace_\eta\cong SO(b_2-3, 2)/SO(2)\times SO(b_2-3).
\]
This is significant, because 
$\Perspace_\eta$ (unlike $\Perspace$) is a
symmetric space. The corresponding result 
for the moduli spaces can be stated
as follows.

\hfill

\corollary\label{_pola_moduli_Corolary_}
Let $(M,\eta)$ be a compact, simple, polarized hyperk\"ahler manifold,
$\bar {\cal M}_{b,\eta}$ a connected component of the weakly polarized 
birational moduli space, defined above,  $G$ the group of integral
orthogonal automorphisms of the lattice $\eta^\bot$ of primitive
elements in $H^2(M)$, and
\begin{equation}\label{_periods_pola_trad_Equation_}
\bar {\cal M}_{b,\eta}\arrow  \Perspace_\eta/G
\end{equation}
the corresponding period map. Then
\eqref{_periods_pola_trad_Equation_}
is a finite quotient. Moreover, 
$\bar {\cal M}_{b,\eta}$ is isomorphic to
a quotient of a symmetric domain
$\Perspace_\eta$ by an arithmetic group
$\Gamma_\eta^I$ acting as above.
\endproof

\hfill

The quotients of such symmetric spaces 
by arithmetic lattices were much studied by
Gritsenko, Hulek, Nikulin, Sankaran and many others
(see e.g. \cite{_GHK:Kodaira_K3_}, \cite{_GHK:moduli_HK_} 
and references therein).
The geometry of $\Perspace_\eta/G$
is in many cases well understood.
Using the theory of automorphic forms,
many sections of pluricanonial (or, in some
cases, plurianticanonical) class can be found,
depending on $q(\eta,\eta)$ and other
properties of the lattice $\eta^\bot$.
In such cases, \ref{_pola_moduli_Corolary_}
can be used to show that the
weakly polarized birational moduli space
has ample (or antiample) canonical class.\footnote{See
\cite{_Debarre_Voisin_} for an alternative
approach to the same problem.}

The automorphic forms on polarized moduli were 
also used to show non-existence of complete
families of polarized K3 surfaces 
(\cite{_BKPS-B:families_}). This program was proposed by 
J. Jorgensen and A. Todorov in 1990-ies, in 
a string of influential (but, sometimes, flawed) preprints,
culminating with \cite{_Jorgenson_Todorov:Borcherds_}.


\section{Hyperk\"ahler manifolds}


In this Section, we recall a number of 
results  about hyperk\"ahler manifolds,
used further on in this paper. For more
details and reference, please see 
\cite{_Besse:Einst_Manifo_}.

\subsection{Hyperk\"ahler structures}
\label{_hk_stru_Subsection_}

\definition
Let $(M,g)$ be a Riemannian manifold, and $I,J,K$
endomorphisms of the tangent bundle $TM$ satisfying the
quaternionic relations
\[
I^2=J^2=K^2=IJK=-\Id_{TM}.
\]
The triple $(I,J,K)$ together with
the metric $g$ is called {\bf a hyperk\"ahler structure}
if $I, J$ and $K$ are integrable and K\"ahler with respect to $g$.

Consider the K\"ahler forms $\omega_I, \omega_J, \omega_K$
on $M$:
\[
\omega_I(\cdot, \cdot):= g(\cdot, I\cdot), \ \
\omega_J(\cdot, \cdot):= g(\cdot, J\cdot), \ \
\omega_K(\cdot, \cdot):= g(\cdot, K\cdot).
\]
An elementary linear-algebraic calculation implies
that the 2-form $\Omega:=\omega_J+\1\omega_K$ is of Hodge type $(2,0)$
on $(M,I)$. This form is clearly closed and
non-degenerate, hence it is a holomorphic symplectic form.

In algebraic geometry, the word ``hyperk\"ahler''
is essentially synonymous with ``holomorphically
symplectic'', due to the following theorem, which is
implied by Yau's solution of Calabi conjecture
(\cite{_Besse:Einst_Manifo_}, \cite{_Beauville_}).

\hfill

\theorem\label{_Calabi-Yau_Theorem_}
Let $M$ be a compact, K\"ahler, holomorphically
symplectic manifold, $\omega$ its K\"ahler form, $\dim_\C M =2n$.
Denote by $\Omega$ the holomorphic symplectic form on $M$.
Suppose that $\int_M \omega^{2n}=\int_M (\Re\Omega)^{2n}$.
Then there exists a unique hyperk\"ahler metric $g$ with the same
K\"ahler class as $\omega$, and a unique hyperk\"ahler structure
$(I,J,K,g)$, with $\omega_J = \Re\Omega$, $\omega_K = \Im\Omega$.
\endproof

\hfill

Further on, we shall speak of ``hyperk\"ahler manifolds''
meaning ``holomorphic symplectic manifolds of K\"ahler
type'', and ``hyperk\"ahler structures'' meaning
the quaternionic triples together with a metric. 

\hfill

Every hyperk\"ahler structure induces a whole 2-dimensional
sphere of complex structures on $M$, as follows. 
Consider a triple $a, b, c\in R$, $a^2 + b^2+ c^2=1$,
and let $L:= aI + bJ +cK$ be the corresponging quaternion. 
Quaternionic relations imply immediately that $L^2=-1$,
hence $L$ is an almost complex structure. 
Since $I, J, K$ are K\"ahler, they are parallel with respect
to the Levi-Civita connection. Therefore, $L$ is also parallel.
Any parallel complex structure is integrable, and K\"ahler.
We call such a complex structure $L= aI + bJ +cK$
{\bf a complex structure induced by a hyperk\"ahler structure}.
There is a 2-dimensional holomorphic family of 
induced complex structures, and the total space
of this family is called {\bf the twistor space}
of a hyperk\"ahler manifold.

\subsection{Bogomolov's decomposition theorem}
\label{_simple_hk_Subsection_}

The modern approach to Bogomolov's decomposition
is based on Calabi-Yau theorem
(\ref{_Calabi-Yau_Theorem_}), Berger's classification
of irreducible holonomy and de Rham's splitting theorem
for holonomy reduction (\cite{_Beauville_},
\cite{_Besse:Einst_Manifo_}). It is worth mentioning
that the original proof of decomposition theorem
(in \cite{_Bogomolov:decompo_}) 
was much more elementary.

\hfill

\theorem\label{_Bogomo_splitting_} 
Let $(M,I,J,K)$ be a compact hyperk\"ahler manifold.
Then there exists a finite covering $\tilde M\arrow M$,
such that $\tilde M$ is decomposed, as a hyperk\"ahler
manifold, into a product
\[
\tilde M = M_1\times M_2 \times \dots M_n \times T,
\]
where $(M_i,I,J,K)$ satisfy $\pi_1(M_i)=0$, $H^{2,0}(M_i,I)=\C$,
and $T$ is a hyperk\"ahler torus. Moreover,
$M_i$ are uniquely determined by $M$ and simply
connected, and $T$ is unique up to isogeny.

\hfill

{\bf Proof:} See \cite{_Beauville_},
\cite{_Besse:Einst_Manifo_}. \endproof

\hfill

\definition
Let $(M,I,J,K)$ be a compact hyperk\"ahler manifold which satisfies 
$\pi_1(M)=0$, $H^{2,0}(M,I)=\C$. Then $M$ is called
{\bf a simple hyperk\"ahler manifold}, or 
{\bf an irreducible hyperk\"ahler manifold}

\hfill

\remark
Notice that \ref{_Bogomo_splitting_} implies that 
irreducible hyperk\"ahler manifolds are simply connected.
In particular, they do not admit a further decomposition.
This explains the term ``irreducible''.

\hfill

As we mentioned in the Introduction, all hyperk\"ahler
manifolds considered further on are assumed to be simple.
Since the Hodge numbers are invariant under K\"ahler deformations,
the deformations of simple manifolds are always simple.
However, the irreducibility is a topological property, as 
implied by the following  lemma.

\hfill

\lemma\label{_irredu_topolo_Lemma_}
Let $M$ be a compact hyperk\"ahler manifold, which is homotopy
equivalent to a simple hyperk\"ahler manifold. Then $M$ is
also simple.

\hfill

{\bf Proof:} 
Let $A^*$ be the part of the rational cohomology of $M$
generated by $H^2(M)$. It is well known 
(see \cite{_Verbitsky:cohomo_} and
\cite{_Verbitsky:coho_announce_})
that $A^*$ is up to the middle dimension a symmetric algebra.
Since $M$ is simply connected, it is 
diffeomorphic to a product of simple hyperk\"ahler
manifolds. Denote by $A^*_i$ the corresponding
subalgebras in cohomology generated by $H^2(M_i)$.
These subalgebras are described in a similar way,
and are symmetric up to the middle.
Then $A^*\cong \bigotimes A_i^*$ by K\"unneth formula. 
Since the algebras $A^*$, $A^*_i$ are symmetric up to the middle, 
this is impossible, as follows from an easy algebraic
computation.
\endproof

\subsection{K\"ahler cone for hyperk\"ahler manifolds}

The following theorem is implied by results of S. Boucksom, 
using the characterization of a K\"ahler cone due
to J.-P. Demailly and M. Paun (see also \cite{_Huybrechts:cone_}). 

Notice that the Beauville-Bogomolov-Fujiki 
form $q$ on \[ H^{1,1}(M,\R):=H^{1,1}(M)\cap H^2(M, \R)\]has signature $(+, -, -,-,...)$, hence the
set of vectors $\nu\in H^{1,1}(M,\R)$
with $q(\nu,\nu)>0$ has two connected components.

\hfill

\theorem\label{_Kahler_cone_Pic=1_Theorem_}
 Let $M$ be a simple hyperk\"ahler manifold
such that all integer $(1,1)$-classes satisfy $q(\nu,\nu)\geq 0$.
Then its K\"ahler cone is one of two components $K_+$ of the set 
$K:=\{ \nu\in H^{1,1}(M,\R)\ \ |\ \ q(\nu,\nu)>0\}$.

\hfill

{\bf Proof:} This is \cite{_V:parabolic_}, 
Corollary 2.6. \endproof

\hfill

For us, the case of trivial Neron-Severi 
lattice is of most interest.

\hfill

\corollary \label{_Kahler_cone_tri_NS_Corollary_}
Let $M$ be a compact, simple hyperk\"ahler manifold
such that $H^{1,1}(M)\cap H^2(M, \Q)=0$. Then
its K\"ahler cone is one of two components  of a set 
$K:=\{ \nu\in H^{1,1}(M,\R)\ \ |\ \ q(\nu,\nu)>0\}.$
\endproof

\subsection{The structure of the period space}
\label{_period_space_grassman_Subsection_}

Let $M$ be a hyperk\"ahler manifold, and $b_2=\dim H^2(M)$.
It is well known that its period space $\Perspace$ 
(see \eqref{_perspace_lines_Equation_})
is diffeomorphic to the Grassmann space 
$Gr_{++}(H^2(M,\R)):=O(b_2-3,3)/SO(2)\times O(b_2-3, 1)$ of
2-dimensional oriented planes $V\subset H^2(M, \R)$
with $q\restrict V$ positive definite.
Indeed, for any line \[ l\in\Perspace\subset {\Bbb P}H^2(M, \C),\]
let $V_l$ be the span of $\langle \Re l, \Im l\rangle$.
From \eqref{_perspace_lines_Equation_} it follows
that $l\cap H^2(M, \R)=0$, hence $V_l$ is an oriented
2-dimensional plane. Since $q(l, \bar l) >0$, the restriction
$q\restrict V$ is positive definite.
This gives a map from $\Perspace$ to $Gr_{++}(H^2(M,\R))$.
To construct the inverse map, we start from
a 2-dimensional plane $V\subset H^2(M, \R)$ and consider the
quadric $\{v\in {\Bbb P}(V\otimes \C) \ \ |\ \ q(v,v)=0\}$.
This quadric is actually a union of 2 points in 
${\Bbb P}(V\otimes \C) \cong\C P^1$, with each of these
points corresponding to a different choice of 
orientation on $V$. This gives an inverse
map from $Gr_{++}(H^2(M,\R))$ to $\Perspace$.

\hfill

The following simple claim is well known. 
For the convenience of the reader, we recall its proof here.

\hfill

\claim \label{_perspace_simply_conne_Claim_}
The period space $\Perspace$ 
is connected and simply connected.

\hfill

{\bf Proof:} We represent $\Perspace$
as $Gr_{++}(H^2(M,\R))=O(b_2-3,3)/SO(2)\times O(b_2-3, 1)$.
The group $O(b_2-3,3)$ is disconnected, but
$O(b_2-3, 1)$ is also disconnected, hence
the connected components cancel each other,
and $Gr_{++}(H^2(M,\R))$ is naturally isomorphic to 
$SO(b_2-3,3)/SO(2)\times SO(b_2-3, 1)$.

To see that it is simply connected, we take a
long exact sequence of homotopy groups
\begin{multline*}
... \arrow \pi_2(Gr_{++}(H^2(M,\R))) 
\arrow 
\pi_1(SO(2)\times SO(b_2-3, 1))\stackrel{(*)}\arrow \\ 
\stackrel{(*)}\arrow\pi_1(SO(b_2-3,3))
\arrow \pi_1(Gr_{++}(H^2(M,\R)))\arrow 0,
\end{multline*}
and notice that the map (*) above is surjective
(it is easy to see from the corresponding maps
of spinor groups and Clifford algebras).
\endproof


\section{Mapping class group of a hyperk\"ahler manifold}


\definition
A connected CW-space $M$ is called {\bf nilpotent}
if its fundamental group $\pi_1(M)$ is nilpotent,
acting nilpotently on homotopy groups of $M$.

\hfill

\definition
Let $M$ be an oriented manifold, 
$\Diff$ the group of oriented diffeomorphisms, and
$\Diff_0$ the group of isotopies, that is,
the connected component of
the group $\Diff$. Then the quotient
$\Diff/\Diff_0$ is called {\bf the mapping class
group of $M$} (see e.g. \cite{_LTYZ:moduli_}).

\hfill 

\definition
 Let $A, A'$ be subgroups in a group $B$.
Recall that $A$ is {\bf commensurable with $A'$}
if $A\cap A'$ has finite index in $A$ and $A'$.
Let $G_\Z$ a group scheme over $\Z$, 
and $G_\R=G_\Q \otimes\Spec\R$
be the corresponding real algebraic group. 
A subgroup $\Gamma \subset G_\R$ is called
{\bf arithmetic} if $\Gamma$ is commensurable
with the group of integer points in $G_\R$.

\hfill

Using rational homotopy theory, formality 
of De\-ligne-Grif\-fiths-Mor\-gan-Sul\-li\-van and Smale's 
h-cobordism, D. Sullivan proved the following general 
result.

\hfill

\theorem\label{_Sullivan_mapping_c_g_Theorem_}
Let $M$ be a compact simply connected (or 
nilpotent)
K\"ahler manifold, $\dim_\C M\geq 3$. Denote by $\Gamma$ the group
of automorphisms of an algebra $H^*(M, \Z)$
preserving the Pontryagin classes $p_i(M)$. 
Then the natural map 
$\Diff/\Diff_0\arrow \Gamma$ has finite kernel,
and its image has finite index in $\Gamma$.
Finally, $\Gamma$ is an arithmetic subgroup in
the group $\Gamma_\Q$ of automorphisms of $H^*(M,\Q)$
preserving $p_i(M)$.

\hfill

{\bf Proof:} 
Theorem 13.3 of \cite{_Sullivan:infinite_}
is stated for general smooth manifolds of $\dim_\R \geq 5$;
to apply it to K\"ahler manifolds, one needs to use
\cite[Theorem 12.1]{_Sullivan:infinite_}. The final
statement is \cite[Theorem 10.3]{_Sullivan:infinite_}.
\endproof

\hfill

For hyperk\"ahler manifolds, the group 
$\Aut(H^*(M, \Q))$ is determined (up to commensurability),
which leads to the following application of
Sullivan's theorem.

\hfill

\theorem\label{_Mapp_cl_group_hk_Theorem_}
Let $M$ be a compact, simple hyperk\"ahler manifold,
its dimension $\dim_\C M=2n$, and $\Gamma_A$ the group
of automorphisms of an algebra $H^*(M, \Z)$
preserving the Pontryagin classes $p_i(M)$. 
Consider the action of $\Gamma_A$ on $H^2(M, \Q)$
and let $\Gamma_2$ be an image of $\Gamma_A$
in $GL(H^2(M, \Q))$. Then 
\begin{description}
\item[(i)] $\Gamma_2$ preserves the Bogomolov-Beauville-Fujiki
form $q$ on $H^2(M, \Q)$.
\item[(ii)] $\Gamma_2$ is an arithmetic subgroup of
$O(H^2(M, \Q), q)$.
\item[(iii)] The natural projection $\Gamma_A\arrow \Gamma_2$
has finite kernel.
\item[(iv)]  The  mapping class
group $\Diff/\Diff_0$ acts on $H^*(M, \Z)$
with finite kernel, and the image 
of $\Diff/\Diff_0$ in $\Gamma_2$ has finite index.
\end{description}
{\bf Proof:} 
From the Fujiki formula 
$v^{2n}= q(v,v)^n$, it is clear that
$\Gamma_A$ preserves the Bogomolov-Beauville-Fujiki,
up to a sign. For $n$ odd, the Fujiki formula fixes the sign,
For $n$ even, the sign is also fixed, because $\Gamma_A$
preserves $p_1(M)$, and (as Fujiki has shown)
$v^{2n-2}\wedge p_1(M)= q(v,v)^{n-1} c$, 
for some $c\in \R$. 
The constant $c$ is positive, because
the degree of $c_2(B)$ is positive
for any Yang-Mills bundle with $c_1(B)=0$
(this argument is based on 
\cite{_Huybrechts:finiteness_}, section 4;
see also \cite{_Nieper:HRR_}).

In \cite[Theorem 13.1]{_Verbitsky:cohomo_} 
(see also \cite[Theorem 2.3]{_Verbitsky:coho_announce_})
it was shown that the group $Spin(H^2(M, \Q), q)$
acts on the cohomology algebra $H^*(M, \Q)$ by
automorphisms, preserving the Pontryagin 
classes.\footnote{In these two papers, the action of
the corresponding Lie algebra was obtained,
giving a $\Spin(H^2(M, \Q), q)$-action by the general
Lie group theory. In \cite[Corollary 8.2]{_V:Mirror_} the action of the centre
of $\Spin(H^2(M, \Q), q)$ was computed. It was shown
that it acts as $-1$ on odd cohomology and trivially
on even cohomology. This Lie algebra, by its construction,
preserves all cohomology classes
which are of type $(p,p)$ for all complex deformations of $M$;
therefore, the group $Spin(H^2(M, \Q), q)$ acts trivially on Pontryagin's
classes.}
Therefore, the group $\Gamma_2\subset O(H^2(M, \Q), q)$ is an
arithmetic subgroup in $O(H^2(M, \Q), q)$.
This gives \ref{_Mapp_cl_group_hk_Theorem_}, (ii).

To see that the kernel $K$ of a map
$\Gamma_A\arrow \Gamma_2$ is finite,
we notice that the subgroup $K\subset \Aut(H^*(M, \Q))$
acts trivially on $H^2(M)$, hence preserves
all Lefschetz ${\goth sl}(2)$-triples
$(L_\omega, \Lambda_\omega, H)$ associated with
different $\omega\in H^{1,1}(M)$. The commutators
of $[L_\omega, \Lambda_\omega]$ generate the Lie algebra
$\goth{so}(H^2(M, \Q), q)$ acting by
derivations on $H^*(M, \Q)$,
as shown in \cite{_Verbitsky:cohomo_} 
(see also \cite{_Verbitsky:coho_announce_}),
hence $K$ centralizes $Spin(H^2(M, \Q), q)$.
The complexification of this group contains
the complex structure operators associated
with any complex, hyperk\"ahler structure on 
$M$ (see \cite{_Verbitsky:cohomo_}, 
\cite{_Verbitsky:coho_announce_}).
Since $K$ centralizes $Spin(H^2(M, \Q), q)$,
$K$ preserves the Hodge decomposition,
for any complex structure $I$ on $M$
of hyperk\"ahler type. Using the Hodge 
decomposition and the Lefschetz 
${\goth sl}(2)$-action, one defines
the Riemann-Hodge pairing, writing down the
Riemann-Hodge formulas as usual;
it is positive definite. Since
$K$ commutes with the ${\goth sl}(2)$-triples
and the Hodge decomposition, it preserves
the Riemann-Hodge pairing $h$. 
Therefore, $K$ is an intersection
of a lattice and a compact group $\Spin(H^*(M), h)$, 
hence finite. We proved 
\ref{_Mapp_cl_group_hk_Theorem_}, (iii).
\ref{_Mapp_cl_group_hk_Theorem_}, (iv)
follows directly from (iii) and
\ref{_Sullivan_mapping_c_g_Theorem_}.
\endproof

\hfill

\remark \label{_Mapping_class_group_lattice_Remark_}
Let $V_\Q$ be a rational vector space equipped
with a quadratic form $q$, and $V_\R:= V_\Q \otimes \R$.
By \cite{_Vinberg_Gorbatsevich_Shvartsman_}, Example 7.5,
the following conditions are equivalent:
\begin{description}
\item[(i)] For any arithmetic subgroup $\Gamma \subset SO(V_\R, q)$,
$\Gamma$ has finite covolume (that is, the quotient
$SO(V_\R, q)/\Gamma$ has finite Haar measure).
\item[(ii)] The algebraic group $SO(V_\Q, q)$ has no
non-trivial homomorphisms to the multiplicative group 
$\Q^{>0}$ of rational numbers (in this case we say that
$SO(V_\Q, q)$ has no non-trivial rational characters).
\end{description}
For $V_\Q=H^2(M, \Q)$ with the Beauville-Bogomolov-Fujiki form,
the latter condition always holds, hence 
the mapping class group is mapped to a discrete 
subgroup of finite covolume
$\Gamma_2\subset SO(H^2(M,\R), q)$.


\section[Weakly Hausdorff manifolds and Hausdorff reduction]
{Weakly Hausdorff manifolds and Hausdorff \\  reduction}


\subsection{Weakly Hausdorff manifolds}

\definition
Let $M$ be a topological space, and $x\in M$ a point.
Suppose that for each $y\neq x$, there exist non-intersecting
open neighbourhoods $U\ni x, V\ni y$. Then $x$ is called
{\bf a Hausdorff point} of $M$.

\hfill

\remark
The topology induced on the
set of all Hausdorff points in $M$ is clearly 
Hausdorff.

\hfill

\definition\label{_weakly_Hau_Definition_}
Let $M$ be an $n$-dimensional 
real analytic manifold, not necessarily Hausdorff. 
Suppose that the set $Z\subset M$ of non-Hausdorff points
is contained in a countable union of real analytic subvarieties
of $\codim\geq 2$. Suppose, moreover, that the
following assumption (called ``assumption {\bf S}''
in the sequel) is satisfied.
\begin{description}
\item[(S)] For every $x\in M$, there is a 
{\em closed} neighbourhood $B\subset M$ of $x$ and a 
continuous map 
$\Psi:\; B \arrow \R^n$ to a closed ball in $\R^n$,
inducing a homeomorphism from an open neighbourhood of $x$
in $B$ onto an open neighbourhood of $\Psi(x)$ in $\R^n$.
\end{description}
Then $M$ is called {\bf a weakly Hausdorff 
manifold}. 

\hfill

\definition
Two points
$x, y\in M$ are {\bf inseparable} (denoted $x\sim y$)
if for any open subsets $U\ni x, V\ni y$, one 
has $U \cap V\neq \emptyset$. 

\hfill

\remark 
A closure of an open
set $U$ contains all points which are inseparable
from some $x\in U$. To extend a homeomorphism
from $\Psi_0:\; B_0 \arrow \R^n$ from an
open neighbourhood $B_0$ to its closure $B$
in order to fulfill the assertion of 
{\bf S} above, we need to extend $\Psi_0$
to all points which are inseparable from
some $x\in B$.

\hfill

\hfill

\remark
Throughout this paper, we could work in much weaker assumptions.
Instead of real analytic, we could demand that $M$ is a 
Lipschitz manifold, and $Z$ has Hausdorff codimension $>1$.
All the proofs in the sequel would remain valid in this general 
situation. Also, the assumption {\bf S} seems to be 
unnecessary, though convenient. In fact, counterexamples
to {\bf S} are hard to find, and it might possibly follow
from the rest of assumptions. 

\hfill

\example\label{_gene_points_in_Teich_Example_}
Let $\Teich$ be a Teichm\"uller space of a hyperk\"ahler
manifold $M$, and $Z\subset M$ the set of all $I\in \Teich$
such that the corresponding Neron-Severi lattice
$H^{1,1}(M,I)\cap H^2(M, \Z)$ has rank $\geq 1$.
Clearly, $Z=\bigcup_\eta Z_\eta$, with the union taken
over all elements 
$\eta \in H^2(M, \Z)$,\footnote{The group $H^2(M, \Z)$
is torsion-free, by the Universal Coefficients Theorem,
because $M$ is simply connected.}
and 
\[ 
  Z_\eta =\{ I\in \Teich \ \  |\ \  \eta \in H^{1,1}(M,I)\}.
\]
As follows from \cite{_Huybrechts:basic_}
(see \ref{_Teich_weakly_Haus_Remark_} below), the complement
$\Teich\!\backslash Z$ is Hausdorff. The period
map $\Per:\; \Teich \arrow \Perspace$ is locally a diffeomorphism,
hence the assumption {\bf S} is also satisfied. Therefore,
$\Teich$ is weakly Hausdorff.

\hfill

The following definition is straightforward;
it is a non-Hausdorff version of a notion
of a manifold with smooth boundary. We 
have to give it in precise detail, because
the notion of a ``boundary'' is ambiguous
in non-Hausdorff situation.

\hfill

\definition\label{_bounda_non_haus_Definition_}
Denote by $\bar U$ the closure of $U$ in $M$, and
by $\bar U^\circ$ the set of interior points of $\bar U$.
Define {\bf the boundary}  as 
$\6_M U:=\bar U \backslash \bar U^\circ$.
We say that an open subset $U\subset M$ of
a smooth manifold $M$ has {\bf smooth boundary,}
if each point in $M$ has a neighbourhood $V$ and a map
mapping the closure of $V$ to $\R^n$, inducing a diffeomorphism
on $V$, and mapping the set 
$\bar V \cap \6_M U $  to $[0, \infty]\times\R^{n-1}$;
The closure $\bar U$ for such $U$ is called {\bf a smooth submanifold
with boundary}. In this case, $\6_M U $ is 
a smooth codimension 1 submanifold of $M$.

\hfill

Further on, we shall need the following claim.
It can be (roughly) stated as follows.
Take a subset $B$ in a weakly Hausdorff $n$-manifold,
diffeomorphic to a closed ball in $U \cong \R^n$
with smooth boundary $\6_UB$. Then its closure $\bar B$ in 
$M$ is obtained by adding two kinds of extra points: 
those in the closure $\overline{\6_U B}$
of $\6_UB$ in $M$ and those which 
are interior to $\bar B$.

\hfill

\claim \label{_closure_of_bou_Claim_}
Let $M$ be a weakly Hausdorff manifold,
$U\subset M$ a subset diffeomorphic to $\R^n$,
and $B \subset U$ a connected, open 
subset of $U$ which has compact closure in $U$ with
smooth boundary $\6_U B\subset U$. Consider the set
$\bar B\backslash \bar B^\circ$ of all points in 
the closure $\bar B$ of $B$ in $M$ which
are not interior in $\bar B$. Then
$\bar B\backslash \bar B^\circ$ 
coincides with the closure $\overline{\6_U B}$ of
$\6_U B$ in $M$. 

\hfill

{\bf Proof:} 
Clearly,  $\6_U B$ contains no interior points of $\bar B$.
Therefore, 
\[
  \overline{\6_U B} \subset \bar B\backslash \bar B^\circ.
\]
We need only to prove the opposite inclusion.

Denote by $W$ the set of Hausdorff points of $M$.
Since $M\backslash W$ has codimension $\geq 2$, 
$W\cap \6_U B$ is dense in $\6_UB$. The boundary $W \cap \6_U B$ separates
$W$ onto two disjoint open subsets, $W_1:= W\cap B$
and $W_2:= W \backslash \bar B$.   
Since $W$ is dense, $\bar B = \bar W_1$, and
$\bar B\backslash \bar B^\circ\subset \bar W_2$.
Therefore, \ref{_closure_of_bou_Claim_} 
would follow if we prove an inclusion
\begin{equation}\label{_bou_in_clo_Equation_}
  \bar W_1 \cap \bar W_2\subset \overline{\6_U B}.
\end{equation}
Let $z\in \bar W_1 \cap \bar W_2$. Then 
in any neighbourhood of $z$ there are
points of $W_1$ and $W_2$.
Since $W$ is a smooth manifold
with countably many codimension $\geq 2$ subvarieties removed,
and $W_1$, $W_2$ are disjoint open subsets of $W$
separated by a smooth boundary $W\cap \6_U B$, this implies
that any neighbourhood of $z$ contains a point in
$W\cap \6_U B$. Indeed, by \ref{_complement_codim_2_conne_Lemma_} 
below, for any path connected open subset $D\subset M$, the
intersection $W\cap D$ is also connected. 
Unless $(W\cap \6 B)\cap D$ is non-empty, the
open set $W\cap D$ is represented as a union
of two non-empty disjoint open subsets $W_1\cap D$
and $W_2\cap D$, which is impossible, because it is connected.
This implies \eqref{_bou_in_clo_Equation_},
and finishes the proof of \ref{_closure_of_bou_Claim_}.
\endproof

\hfill

The following trivial lemma, used in the proof of
\ref{_closure_of_bou_Claim_}, is well-known; we include it
here for completeness.

\hfill

\lemma\label{_complement_codim_2_conne_Lemma_}
Let $M$ be a path connected real analytic manifold, and 
$W= M \backslash \bigcup Z_i$, where $\bigcup Z_i$
is a union of countably many real analytic manifolds
of codimension at least 2. Then $W$ is path connected.

\hfill

{\bf Proof:} This result is clearly local. Therefore,
we may assume that $M$ is isomorphic to $\R^n$. Given two 
points $x, y \in W$, we shall prove that there is $z\in W$
such that a straight segment of a line 
connecting $z$ to $x$ and the one connecting 
$z$ to $y$ belong to $W$. Let $P^x\cong {\R P^n}$ be the set
of all lines passing through $x$, and $P_W^x$ the set
of these lines which belong to $W$. Clearly,
the set $P_{Z_i}^x$ of lines $l\in P$ intersecting $Z_i$,
being a projection of $Z_i$ to $P$, 
has real codimension 1 in $P$. Therefore, the complement to a 
set $P_W^x$ is of measure 0 in $P^x$. Similarly one defines
$P_W^y$ and proves that it is dense. Let now
$Q$ be the set of all pairs of 
lines $l^x \in P^x, l^y\in P^y$ which intersect.
Clearly, $Q$ is equipped with smooth projections
$\pi_x$, $\pi_y$ to $P^x$ and $P^y$,
with 1-dimensional fibers. Since
the complements to $P_W^x$ and $P_W^y$ in $P^x$ and $P^y$
have measure 0, the intersection $\pi_x^{-1}(P_W^x) \cap \pi_y^{-1}(P_W^y)$
is non-empty. For each pair of lines
\[
  (l^x, l^y)\in \pi_x^{-1}(P_W^x) \cap \pi_y^{-1}(P_W^y)\subset Q,
\]
$l^x$ and $l^y$ are lines which are contained in $W$,
intersect and connect $x$ to $y$.
\endproof

\subsection{Inseparable points in weakly Hausdorff manifolds}

\lemma\label{_interse_interior_Lemma_}
Let $M$ be a weakly Hausdorff manifold, $x, y \in M$
inseparable points, and $U\ni x, V\ni y$ open sets.
Then $x$ and $y$ are interior points of 
$\bar U \cap \bar V$, where $\bar U, \bar V$
denotes the closure of $U$, $V$. 

\hfill

\remark
This statement is false without the weak Hausdorff assumption.
Indeed, take as $M$ the union of two real 
lines, with $t<0$ identified, $x$ the 0 of the first line,
$y$ the 0 of the second line. Choose a neighbourhood $U$ 
of $x$ and $V$ of $y$. The points $x$ and $y$ are clearly
inseparable, but the intersection of $\bar U \cap \bar V$
is a closure of an interval $[-a, 0[$, with $a>0$ a positive number,
hence $x$ and $y$ are not interior points of $\bar U \cap \bar V$.

\hfill

{\bf Proof of \ref{_interse_interior_Lemma_}:}
Consider an open ball $B\subset U$ with smooth boundary $\6_U B$
containing $x$. Since  $x$ and $y$ are 
inseparable, $y$ belongs to a closure $\bar B$ of $B$. 
Then either $y$ is interior
in $\bar B$, or $y$ lies in the closure
of its boundary $\6_U B$, as follows from \ref{_closure_of_bou_Claim_}.
To prove \ref{_interse_interior_Lemma_}
it remains to show that the second option
is impossible. Using the assumption ``{\bf S}''
of the definition of weakly Hausdorff manifolds,
we obtain that $\Psi(y)=\Psi(x)$,
where $\Psi:\; B \arrow \R^n$
is the map defined in ``{\bf S}''.
Choosing $B$ sufficiently small,
we can always assume that $\Psi\restrict B$ 
is a homeomorphism. Then $\Psi(x)=\Psi(y)$ is
in the interior of $\Psi(B)$, hence
$\Psi(y)\notin \Psi(\6_U B)$. Since
$\Psi$ is continuous, $\Psi^{-1}(\Psi(\6_U B))$
contains the closure of $\6_U B$. Therefore,
$y \notin \overline{\6_U B}$ by 
\ref{_closure_of_bou_Claim_}.
We proved \ref{_interse_interior_Lemma_}. 
\endproof

\hfill

We shall also need the following trivial lemma.

\hfill

\lemma\label{_interse_with_Hausdo_Lemma_}
In assumptions of \ref{_closure_of_bou_Claim_}, let
$W\subset M$ the set of Hausdorff points of $M$.
Then the intersection $W \cap \bar B^\circ$ lies in $B$.

\hfill

{\bf Proof:} Denote by $B_{cl}$ the union of $B$ and $\6_U B$
(\ref{_bounda_non_haus_Definition_}).
Let $x\in W\cap \bar B^\circ$.
Then $x$ is a limit of a sequence $\{x_i\} \in B$.
Since $B_{cl}$
is compact, $\{x_i\}$ has a limit point $x'$ in $B_{cl}$.
Since $x$ is Hausdorff, and $x\sim x'$, one has $x=x'$.
Therefore, $x\in B_{cl}$.
By \ref{_closure_of_bou_Claim_}, one and
only one of two things happens:
either $x$ is interior in  $\bar B$,
or it belongs to the closure $\overline{\6_U B}$ of 
the smooth sphere $\6_U B= B_{cl}\backslash B$.
The later case is impossible, because $x$ is 
interior in $\bar B$. Therefore, $x$ is interior in $B_{cl}$.
\endproof

\hfill

\proposition\label{_insepara_equiva_Proposition_}
Let $M$ be a weakly Hausdorff manifold, and $\sim$ be 
inseparability relation defined above. Then $\sim$ is
an equivalence relation.

\hfill

\remark
Without the weak Hausdorff assumption, 
$\sim$ is not an equivalence relation. Indeed,
consider for example a union $\R \coprod \R \coprod \R$
of three real lines and glue $t<0$ for the first two lines, and $t>0$ for
the second two. Then $0_1$ (the zero on the first line)
is inseparable from $0_2$, and $0_2$ from $0_3$, 
but $0_1 \not\sim 0_3$.

\hfill

{\bf Proof of \ref{_insepara_equiva_Proposition_}:} 
Only transitivity needs to be proven. 
Let $x_1\sim x_2$, $x_2 \sim x_3$ be points in $M$,
$U_1 \ni x_1$, $U_3 \ni x_3$ their neighbourhoods.
By \ref{_interse_interior_Lemma_},
$x_2$ is an interior point of $\bar U_1$ and $\bar U_3$.
Therefore, $\bar U_1 \cap \bar U_3$ is non-empty,
and contains an open subset $A:=\bar U_1^\circ \cap \bar U_3^\circ$,
where $\bar U_i^\circ$ be the set of interior points
of $\bar U_i$. The intersection
$A \cap W$ of $A$ with the set of Hausdorff points
in non-empty, because $W$ is dense. The intersection 
$\bar U_i^\circ\cap W$ lies in $U_i$,
as follows from \ref{_interse_with_Hausdo_Lemma_}, 
hence $A\cap W$ lies in $U_1$ and $U_3$, and these two
open sets have non-trivial intersection.
\endproof

\hfill 

Further on, we shall be interested in 
the quotient $M/\!{}_\sim$, equipped with a quotient topology.
By definition, a subset  $U\subset M/\!{}_\sim$ is open
if its preimage in $M$ is open, and closed if its
preimage in $M$ is closed. 

\hfill

\claim\label{_open_insepara_Claim_}
Let  $M$ be a weakly Hausdorff manifold, 
and $B\subset M$ an open subset with smooth
boundary. Consider its closure $\bar B$, and let 
$\bar B^\circ$ be the set of its interior points.
Then $\bar B^\circ$ is the set of all points
$y\in M$ which are inseparable from some
$x\in B$.

\hfill

{\bf Proof:}
Let $x\in  B$ be any point,
and $y \in M$ a point inseparable from $x$.
By \ref{_interse_interior_Lemma_}, 
for any neighbourhood $U\ni y$, 
$y$ is an interior point of $\bar U\cap \bar B$.
Therefore, $y$ is an interior point of $\bar B$.

To finish the proof of \ref{_open_insepara_Claim_},
it remains to show that any interior point  $z\in \bar B$
is inseparable from some $z'\in B$. This statement is local
in $B$, hence we may assume that $B$ is diffeomorphic to an
open ball satisfying the assumptions of the property
{\bf S} of \ref{_weakly_Hau_Definition_}.

 Choose a diffeomorphism
$B \stackrel \Psi \arrow B^{\R^n}$ to an open ball
$B^{\R^n} \subset \R^n$. Using the property
{\bf S} of \ref{_weakly_Hau_Definition_},
we may assume that $\Psi$ can be extended
to a continuous map from the closure
$\bar B$ to the closed ball $\bar B^{\R^n}$.

Any interior point  $z\in \bar B$
can be obtained as a limit of a sequence
of points  $\{z_i\}\subset B$. 
Let $\zeta\in \bar B^{\R^n}$ be a 
limit of $\{\Psi(z_i)\}$ in $\bar B^{\R^n}$,
which exists because $\bar B^{\R^n}$ is compact.
Choosing a subsequence, we may also assume that
$\lim\{\Psi(z_i)\}$ is unique.
Then $\zeta =\Psi(z)$, and it is an
interior point of $\bar B^{\R^n}$, as follows
from \ref{_closure_of_bou_Claim_}.
Since $B \stackrel \Psi \arrow B^{\R^n}$ 
is a diffeomorphism, the sequence
$\{z_i\}$ has a limit $z'\in B$.
Since $\Psi(z)=\Psi(z')=\lim\{\Psi(z_i)\}$, 
the point $z$ is inseparable from $z'$.
\endproof

\hfill

\theorem\label{_Hausdorff_redu_Theorem_}
Let $M$ be a weakly Hausdorff manifold, and
$\sim$ the inseparability relation.
Consider the quotient space $M/\!{}_\sim$ equipped
with a natural quotient topology. Then 
$M/\!{}_\sim$ is Hausdorff, and the projection
map $M \stackrel \phi \arrow M/\!{}_\sim$ is open.

\hfill

{\bf Proof:} Since $M$ is a manifold, we can choose 
a base of open subsets $U\subset M$ with smooth boundary. 
By \ref{_open_insepara_Claim_}, $\phi^{-1}(\phi(U))= \bar U^\circ$,
where $\bar U^\circ$ is the set of all interior points
of the closure $\bar U$. Therefore, the image of $U$ is 
open in $M/\!{}_\sim$, and $\phi$ is an open map. 

Denote by $\Gamma_\sim\subset M\times M$ the graph of $\sim$.
It is well known that a topological space $X$ is Hausdorff
if and only if the diagonal $\Delta$ is closed in $X\times X$.
Since the projection 
$M\times M \stackrel {\phi\times \phi}\arrow  
M/\!{}_\sim \times M/\!{}_\sim$
is open, and 
\[ \phi(M\times M\backslash \Gamma_\sim)= 
   \left(M/\!{}_\sim \times M/\!{}_\sim\right)\backslash \Delta,
\]
to prove that $M/\!{}_\sim$ is Hausdorff it remains to show that
$\Gamma_\sim$ is closed in $M\times M$.

Let $(x,y)\notin \Gamma_\sim$, equivalently, $x\not\sim y$.
Choose open neighbourhoods $U\ni x, V\ni y$,
$U\cap V=\emptyset$. Then $U \times V \cap \Gamma_\sim=\emptyset$.
This implies that $\Gamma_\sim$ is closed. We proved that
$M/\!{}_\sim$ is Hausdorff. \endproof

\subsection{Hausdorff reduction for weakly Hausdorff manifolds}
\label{_Hausdorff_red_Subsection_}

\definition
Let $X \stackrel \phi \arrow Y$ be a surjective morphism of
topological spaces, with $Y$ Hausdorff. Suppose that for
any map $X \stackrel {\phi'}\arrow Y'$, with $Y'$ Hausdorff, the map
$\phi'$ is factorized through $\phi$.
Then $\phi$ is called {\bf the Hausdorff reduction map},
and $Y$ {\bf the Hausdorff reduction of $X$}.
Being an initial object in the category of diagrams
$X\stackrel {\phi'}\arrow Y'$ (with $Y'$ Hausdorff),
the Hausdorff reduction if obviously unique, if it exists.

\hfill

\remark
If $x\sim y$ are inseparable points of $M$, any morphism
$M\stackrel \phi \arrow M'$ to a Hausdorff space $M'$
satisfies $\phi(x)=\phi(y)$. Therefore, whenever 
the quotient $M/\!{}_\sim$ is Hausdorff, it is a Hausdorff
reduction of $M$.

\hfill

\example 
By \ref{_Hausdorff_redu_Theorem_}, for any weakly Hausdorff
manifold $M$, the quotient $M/\!{}_\sim$ is its Hausdorff reduction.

\hfill

\definition
A {\bf local homeomorphism} is
a continuous map $X\stackrel \psi\arrow Y$
such that for all $x\in X$ there is a neighbourhood
$U\ni x$ such that $\psi\restrict U$ is a homeomorphism
onto its image, which is open in Y. If $\psi$ is
also a smooth, it is called {\bf a local diffeomorphism},
or {\bf etale map}.

\hfill

\theorem\label{_Teich_b_manifold_Theorem_}
Let $M$ be a weakly Hausdorff manifold, and
\[ \phi:\; M \arrow M/\!{}_\sim\] its Hausdorff reduction.
Then $\phi$ is etale,
and $M/\!{}_\sim$ is a Hausdorff manifold.

\hfill

{\bf Proof:} Let $U \subset M$ be an open neighbourhood
of a given point $x$, diffeomorphic to $\R^n$,
and $B\subset U$ a closed neighbourhood diffeomorphic
to a closed ball. Since $U$ is Hausdorff, the
restriction $\phi\restrict U$ is injective. 
An injective map from a compact $B$ to a Hausdorff
space is a homeomorphism to its image. Then
the restriction of $\phi$ to interior of $B$
is a homeomorphism.  \endproof

\subsection{Inseparable points in the marked moduli and the 
Teichm\"uller space}

\definition
Let $K\subset \Gamma$ be a subgroup of the mapping class group
acting trivially on $H^2(M)$. It is a finite group by 
\ref{_Mapp_cl_group_hk_Theorem_} (iv).
Recall that {\bf the marked moduli space} $\Teich_{H^2}$ is a quotient
of the Teichm\"uller space by $K$.

\hfill

The following result is due to D. Huybrechts.

\hfill

\theorem\label{_bira_points_inseparable_Theorem_}
(\cite{_Huybrechts:cone_})
Let $M$ be a hyperk\"ahler manifold, $\Teich_{H^2}$ its marked moduli 
space, and $x, y \in \Teich_{H^2}$  points corresponding to
hyperk\"ahler manifolds $M_x$ and $M_y$.
Suppose that $x$ and $y$ are inseparable, in the sense of 
\ref{_insepara_Definition_}. Then
the manifolds $M_x$ and $M_y$ are bimeromorphically
equivalent. Conversely, if $M_1$ and $M_2$
are bimeromorphically equivalent, they can
be realised as inseparable points on the
Teichm\"uller space. \endproof

\hfill

We extend this result to a Teichm\"uller space.

\hfill

\proposition\label{_non_Hausdorff_marked_and_Teich_Proposition_}
Let $I\in \Teich$ be a non-Hausdorff point. Then
its image in $\Teich_{H^2}$ is also non-Hausdorff.

\hfill

\ref{_non_Hausdorff_marked_and_Teich_Proposition_}
is proven at the end of this section. 
Its proof easily follows from

\hfill

\theorem\label{_stabili_in_H^2_maping_class_Theorem_}
Let $I\in \Teich$, and $K_I\subset K$ be a stabilizer of $I$
in the group $K$, defined as a subgroup of the mapping class group
acting trivially on $H^2(M)$. Then 
\begin{description}
\item[(i)]
Denote by $G_I\subset K$ be a subgroup fixing a connected component
$\Teich^I\ni I$ of $\Teich$. Then $K_I$ is normal in $G_I$.  

\item[(ii)] The group $K_I$ is independent
from the choice of $I$ in the component $\Teich^I$
of the Teichm\"uller space. 

\item[(iii)]
Moreover, the natural projection
$\Teich^I\arrow \Teich_{H^2}^I$ is a finite covering, with
$\Teich_{H^2}^I=\Teich^I/(G_I/K_I)$. Here, $\Teich_{H^2}^I$ denotes the
connected component of $\Teich_{H^2}$ containing $I$.
\end{description}

{\bf Proof:} Part (i) clearly follows from (ii). 
To prove (ii), consider the action of $K_I$ on $T_I\Teich$.
We shall prove that a stabilizer $St(K_I)\subset \Teich^I$
is open and closed in $\Teich^I$.

By Bogomolov's theorem (\ref{_bogomo_smooth_Theorem_}),
$T_I\Teich$ is naturally identified with $H^{1,1}_I(M)$.
However, the action of $K$ on $H^2(M)$ is trivial, hence
any $\alpha\in K_I$ acts trivially on $T_I\Teich$. For
any finite order diffeomorphism of $\Teich$, dimension
of its fixed point set passing through $I$ is equal
to the dimension of the corresponding unit eigenspace
in $T_I\Teich$. Therefore, $K_I$ acts as identity
on an open subset $St(K_I)\subset \Teich^I$. 
This subset is also closed, which can be seen from
the following argument. 

Let $J\in \Teich^I$ be a point,
and $\{J_k\}\in St(K_I)\subset \Teich^I$ a sequence converging
to $J\in \Teich^I$. To prove that $St(K_I)$ is closed,
it would suffice to show that $J\in St(K_I)$.

Consider some $h\in K_I$, and
let $\tilde h$ be the lift of $h$ to $\Diff(M)$. 
To prove that $J\in St(K_I)$, it would suffice to find
an isotopy $\nu$ such that $\nu^* J= \tilde h ^*J$.

Consider at infity-dimensional space of integrable almost 
complex structures $\Comp$, equipped with topology of Fr\'echet
convergence (it is a Fr\'echet manifold:
see \cite{_Hamilton:Nash_}). The space $\Teich$ is
defined as a quotient $\Teich:=\Comp/\Diff^0$,
with factor topology (\cite{_Hamilton:Nash_}, 
\cite{_Kuranishi:new_}). 

Chose representatives $\tilde J$, $\tilde J_k\in \Comp$.
From the definition of the topology on $\Teich$
it is obvious that we can chose $\{\tilde J_k\}$
converging to $\tilde J$ in the Fr\'echet topology
on $\Comp$.

Choose a sequence  $[\omega_k]\in H^2(M, J_k)$ of 
K\"ahler classes converging to a K\"ahler class
$[\omega]$ on $(M,J)$. This is possible to do by Kodaira's
deformational stability of K\"ahler structures.
Let $\omega_k, \omega$ be the corresponding
Ricci-flat K\"ahler metrics on $(M, \tilde J_k)$ and
$(M, \tilde J)$.
Since the solutions of Calabi-Yau problem continuously
depend on the data (complex structure and the K\"ahler class),
the tensors $\omega_k$ converge to $\omega$.

Chose isotopies $\nu_k\in \Diff_0$ in such a way that
$\nu_k^* \tilde J_k = \tilde h^* \tilde J_k$.
This is possible to do, because $\tilde h^* \tilde J_k$
is equivalent to $\tilde J_k$ in the Teichm\"uller space.
The map $(\tilde h^{-1} \nu_k)^*$ is holomorphic
on $(M,\tilde J_k)$ and preserves the complex structure,
hence induces an isometry of Calabi-Yau metrics.
Since the group of isometries is compact, 
$\tilde h^{-1} \nu_k$ converges to a holomorphic
isometry $N$ of $(M,\tilde J)$. Let $\nu:= \tilde h N$.
By construction, $\nu = \lim \nu_k$ in Fr\'echet topology.
Since $\Diff^0$ is closed in $\Diff$, and all $\nu_k$
lie in $\Diff^0$, this implies that $\nu \in \Diff^0$.
We have just shown that $\nu^* J= \tilde h ^*J$,
where $\nu$ is an isotopy.

Now, $\Teich^I$ is a union of subsets $St(K_{I'})$, for
various $I'\in \Teich^I$, which are all open and closed in
$\Teich^I$. Since $\Teich^I$ is connected, this implies
$\Teich^I = St(K_I)$.
 This proves \ref{_stabili_in_H^2_maping_class_Theorem_}
(ii).

To prove \ref{_stabili_in_H^2_maping_class_Theorem_}
(iii), recall that $K$ is defined as a kernel of the natural map
$\Gamma\arrow GL(H^2(M, \R))$, where $\Gamma:= \Diff^+(M)/\Diff_0(M)$
is the mapping class group. From \ref{_stabili_in_H^2_maping_class_Theorem_}
(ii) we obtain that $\Teich^I\arrow \Teich_{H^2}^I$ is a quotient of $\Teich^I$
by $G_I\subset K$. Then \ref{_stabili_in_H^2_maping_class_Theorem_}
(iii) is implied by the following trivial lemma, applied to $X=\Teich^I$
and $G=G_I/K_I$.

\hfill

\lemma
Let $G$ be a finite group freely acting on a manifold
$X$ (possibly non-Hausdorff). Assume that $X/G$ is also
a manifold, and the projection $X \arrow X/G$ is locally
a homeomorphism. Then the quotient map $X \stackrel \pi \arrow X/G$
is a covering.

\hfill

{\bf Proof:} It would suffice to prove that $G$ acts properly,
that is, to show that each point $x\in X$ has a neighbourhood $U\ni x$ 
which is disjoint from $gU$, for all non-trivial $g\in G$.

Since the projection map $X \stackrel \pi \arrow X/G$ is etale, there
exists a neighbourhood $U\ni x$ such that $\pi:\; X \arrow X/G$ is
a homeomorphism. The group $G$ freely acts on $\bigcup_{g\in G} gU$,
and for each $g\in G$ the restriction $\pi:\; gU \arrow \pi(U)$
is a homeomorphism. Then
the sets $gU$ never intersect, for different $g\in G$.
Therefore, $U$ is a neighbourhood of $x$ which satisfies
$\forall g\in G, g\neq e, U\cap gU = \emptyset.$ 
\endproof

\hfill

{\bf Proof of \ref{_non_Hausdorff_marked_and_Teich_Proposition_}:}
\ref{_non_Hausdorff_marked_and_Teich_Proposition_}
follows from \ref{_stabili_in_H^2_maping_class_Theorem_},
because $G_I/K_I$ is a finite group which properly acts
on $\Teich^I$, hence
the projection $\Teich^I\arrow\Teich_{H^2}^I$ is a finite covering,
which maps non-Hausdorff points to non-Hausdorff points.
\endproof

\subsection{The birational Teichm\"uller space for a hyperk\"ahler manifold}

\remark\label{_Teich_weakly_Haus_Remark_}
Let $M_1,M_2$ be bimeromorphically equivalent 
hyperk\"ahler manifolds. By \cite[Proposition 9.2]{_Huybrechts:basic_}
and \ref{_non_Hausdorff_marked_and_Teich_Proposition_}, 
the Neron-Severi lattice
$\NS(M_i) = H^{1,1}(M,\Z)$ 
has rank $\geq 1$, unless the bimeromorphism
$M_1 \rightsquigarrow M_2$ is biregular. Therefore, a point
$I \in \Teich$ with $\rk \NS(M,I)=0$
must be separable. This argument
was used earlier in \ref{_gene_points_in_Teich_Example_}
to prove that $\Teich$ is weakly 
Hausdorff.

\hfill

Clearly, the map $\Per:\; \Teich_b\arrow \Perspace$
is well defined (it follows directly from the
definition of the Hausdorff reduction). Indeed, the birational
Teichm\"uller space $\Teich_b$ is obtained as a Hausdorff reduction
of the Teichm\"uller space.
The main result of this paper is the following theorem

\hfill

\theorem\label{_global_Torelli_Theorem_}
(global Torelli theorem)
Let $M$ be a simple hyperk\"ahler manifold, 
$\Teich_b$ its birational Teichm\"uller space,
and 
\begin{equation}\label{_bira_peri_map_Equation_}
\Per:\; \Teich_b\arrow \Perspace
\end{equation} 
the period
map defined as above. Then
\eqref{_bira_peri_map_Equation_} is a
diffeomorphism, for each connected 
component of $\Teich_b$.

\hfill

\ref{_global_Torelli_Theorem_} follows from
\ref{_Per_covering_Proposition_}, because
$\Perspace$ is simply connected (\ref{_perspace_simply_conne_Claim_}).

\hfill

\proposition\label{_Per_covering_Proposition_}
Consider the map $\Per:\; \Teich_b\arrow \Perspace$
defined as in \ref{_global_Torelli_Theorem_}.
Then $\Per$ is a covering.

\hfill

\ref{_Per_covering_Proposition_}
is proven in \ref{_Per_cove_Proposition_proof_Remark_} below.

\hfill

As an immediate corollary, we obtain the following result

\hfill

\corollary
Let $M$ be a hyperk\"ahler manifold, $\Teich_{H^2}$ its marked moduli 
space, and $\Teich\stackrel \Psi\arrow \Teich_{H^2}$ the natural projection.
Then $\Psi$ is a diffeomorphism on each connected component.

\hfill

{\bf Proof:} By \ref{_stabili_in_H^2_maping_class_Theorem_},
$\Psi$ is a finite covering, hence it induces a finite covering
of the corresponding Hausdorff reductions. However,
$\Psi$ induces an isomorphism of Hausdorff reductions,
because each component of $\Teich_{H^2}/\sim$ and 
 $\Teich/\sim$ is isomorphic to $\Perspace$. \endproof

\hfill

\remark
The Hausdorff reduction $\Teich/\!{}_{\sim}$ classifies
complex structures on $M$ up to ``bimeromorphic
equivalence'' and the action of the isotopy group.
We call $\Teich/\!{}_\sim$ {\bf the birational
Teichm\"uller space}, denoting it as $\Teich_b$.
However, the term ``bimeromorphic equivalence'' is vague.
Clearly, there are distinct points in $\Teich/\!{}_\sim$
which represent bimeromorphic (and biholomorphic)
hyperk\"ahler manifolds. A better description of
this equivalence might be gleaned from
\cite{_Huybrechts:cone_} and \cite{_Boucksom_} 
(I am grateful to Eyal Markman for this observation,
also found in \cite{_Markman:survey_}).
Consider the Hodge isometry 
$f:\; H^2(M_1, \Z) \arrow  H^2(M_2, \Z)$ between the
second cohomology corresponding to two inseparable
points in $\Teich$. In the language of Boucksom, $f$
maps the K\"ahler cone to one of the ``rational chambers''
of the positive cone. As shown in \cite[Theorem 4.3]{_Boucksom_}, 
there are three possibilities:
\begin{description}
\item[(i)] $f$ could map the K\"ahler cone to the K\"ahler
cone, which means that $f$ is induced by an isomorphism. 
In this case $M_1$ and $M_2$ correspond to the same
points of the marked moduli space.
\item[(ii)] $f$ could map the K\"ahler cone onto a different
rational chamber, which belongs to the {\em fundamental 
uniruled chamber}. In this case $f$ is induced by a graph 
of a bimeromorphic morphism.
\item[(iii)] $f$ could map the K\"ahler cone onto a different
rational chamber, which does {\bf not} belong to the 
fundamental uniruled chamber. In this case $f$ is induced
by a {\bf reducible} correspondence. One of its irreducible
components is a graph of a birational morphism. Other
components, which necessarily exist, will appear as
certain fiber products of uniruled divisors. 
\end{description}


\section{Hyperk\"ahler lines in the moduli space}


In this Section, we introduce generic hyperk\"ahler
lines (GHK lines), and prove that every two
points of the Teichm\"uller space of hyperk\"ahler manifold
are connected by a sequence of 5 GHK lines
(\ref{_4_GHK_lines_Proposition_})

\subsection{Hyperk\"ahler lines and hyperk\"ahler structures}

\definition
Let $M$ be a simple hyperk\"ahler manifold, $\Perspace$
its period space, and $W\subset H^2(M, \R)$ an oriented  
3-dimensional subspace, such that $q\restrict W$ is positive definite.
Consider a 2-dimensional sphere $S_W\subset \Perspace$
consisting of all oriented 2-dimensional planes $V\subset W$.
Using an isomorphism \[ \Perspace\cong Gr_{+,+}(H^2(M, \R))\] constructed in
Subsection \ref{_period_space_grassman_Subsection_}, 
we can consider $S_W$ as a subvariety in $\Perspace$.
This subvariety is called
{\bf a hyperk\"ahler line associated with a 3-dimensional
plane $W\subset H^2(M, \R)$}.

\hfill

\remark\label{_sphere_hk_Remark_}
Let $(M,g,I,J,K)$ be a hyperk\"ahler structure,
$S\subset \Teich$ the sphere of induced complex structures
defined as in Subsection \ref{_hk_stru_Subsection_}, 
and $W:=\langle \omega_I, \omega_J, \omega_K\rangle 
\subset H^2(M, \R)$ the corresponding 3-dimensional plane.
It is easy to see that the sphere $\Per (S) \subset \Perspace$
coincides with the hyperk\"ahler line $S_W$ defined as above.
This explains the term.

\hfill

\definition\label{_GHK_line_Definition_}
Let $S_W\subset \Perspace$ be a hyperk\"ahler line
associated with a 3-di\-men\-si\-onal subspace $W\subset H^2(M, \R)$. 
We say that $S_W$ is {\bf a generic hyperk\"ahler line}
if the orthogonal complement to $W$ has no rational points: 
\[ W^\bot \cap H^2(M, \Q)=0.\]
Often, we shall abbreviate ``generic hyperk\"ahler line'' to
``GHK line''

\subsection{Generic hyperk\"ahler lines and the Teichm\"uller space}

Let $(M,I)$ be a hyperk\"ahler manifold.
The Hodge structure on $H^2(M, I)$ is determined from the
Bogomolov-Beauville-Fujiki form $q$ and the
corresponding 1-dimensional space $l=\Per(I)\subset H^2(M, \C)$:
one has $H^{2,0}(M,I)=l$, $H^{0,2}(M,I)=\bar l$, and
$H^{1,1}(M,I)= \langle l, \bar l\rangle^\bot$,
where $\bot$ denotes the orthogonal complement.
We define the Neron-Severi lattice of $(M,I)$
as $\NS(M,I):=H^{1,1}(M,I) \cap H^2(M, \Z)$. 
Since $H^{1,1}(M,I)= \langle l, \bar l\rangle^\bot$,
the lattice $\NS(M,I)$ depends only on the
point $\Per(I)\in \Perspace$. We shall often consider
the Neron-Severi lattice of a point $l\in \Perspace$,
defined as above.  Since a simple hyperk\"ahler manifold
is simply connected, $\NS(M,I)=\Pic(M,I)$. This allows us to define
the Picard group $\Pic(l)$ for $l\in \Perspace$:
\[
  \Pic(l)=NS(l)=\langle l, \bar l\rangle^\bot\cap H^2(M, \Z).
\]

\hfill

\claim\label{_GHK_lines_and_Pic_Claim_}
Let $S\subset \Perspace$ be a hyperk\"ahler line,
associated with a 3-\-di\-men\-sio\-nal subspace 
$W\subset H^2(M, \R)$. 
Then the following assumptions are equivalent.
\begin{description}
\item[(i)] $S$ is a GHK line
\item[(ii)] For some $l\in S$, the corresponding 
Neron-Severi lattice $\NS(M,l)$ is trivial.
\item[(iii)] For some $w\in W$, its orthogonal complement 
$w^\bot\subset H^2(M, \R)$ has no non-zero rational points.
\end{description}
{\bf Proof:} The points of $S$ are parametrized by 
oriented 2-dimensional planes $V\subset W$, and
the corresponding Neron-Severi lattice $\NS(M,V)$
is $V^\bot \cap H^2(M, \Z)$. Now, the chain of inclusions
\[
 W^\bot \cap H^2(M, \Q)
\subset V^\bot \cap H^2(M, \Q)\subset w^\bot\cap H^2(M, \Q)
\]
immediately brings the implications
(iii) $\Rightarrow$ (ii) $\Rightarrow$ (i).
To finish the proof, it remains to deduce (iii) from
(i). Let 
\[ R:= 
\bigcup\limits_{\stackrel{\eta\in H^2(M, \Q)}{\eta\neq 0}}\eta^\bot
\]
be the union of all hyperplanes orthogonal to non-zero
rational vectors. Since $W^\bot \cap H^2(M, \Q)=0$,
$W$ does not lie in $R$. Therefore, $W\cap R$
is a countable union of planes of positive
codimension. Take $w\in W \backslash R$. Clearly,
$w^\bot\cap H^2(M, \Q)=0$.
\endproof

\hfill

\remark\label{_gene_on_hk_line_Remark_}
The same proof also implies that for any generic 
hyperk\"ahler line, the set of all $I\in S$
with $\NS(M,I)\neq 0$ is countable. Indeed, it
is a countable union of closed complex subvarieties
of positive codimension in $\C P^1$.

\subsection{Connected sequences of GHK lines}

Further on in this subsection, we shall use the following
trivial linear-al\-ge\-bra\-ic lemma.

\hfill

\lemma\label{_vector_posi_def_Lemma_}
Let $A$ be a real vector space, equipped with
a non-degenerate scalar product $q$, 
$W\subset A$ a $d$-dimensional subspace on which $q$ 
is positive definite,\footnote{Further
on, such spaces will be called {\bf positive}.}
and $W'\subset A$ a positive subspace of dimension
$d'<d$. Then there exists a non-zero vector 
$b\in W$, such that the subspace generated 
by $b$ and $W'$ is also
positive.

\hfill

{\bf Proof:} Assume that $W\cap W' =0$
(otherwise, we could choose $b\in W\cap W'$).
Then $\dim {W'}^\bot \cap W>0$. Choose
$b\in {W'}^\bot \cap W$. \endproof

\hfill

\remark
Of course, the set of such $b$ is open in $W$.

\hfill

\proposition\label{_4_GHK_lines_Proposition_}
Let $x, y \in \Perspace$. Then $x$ can be connected
to $y$ by a sequence of 5 sequentially intersecting
GHK lines.

\hfill

{\bf Proof. Step 0:} Using the identification between
$\Perspace$ and the Grassmann space $Gr_{+,+}(H^2(M, \R))$
(Subsection \ref{_period_space_grassman_Subsection_}),
we shall consider points of $\Perspace$ as 
2-dimensional subspaces  $V\subset H^2(M, \R)$
with $q\restrict V$ positive definite. The hyperk\"ahler
lines are understood as 3-dimensional spaces
$W\subset H^2(M, \R)$
with $q\restrict W$ positive definite.
Under this identification, the incidence relation 
is translated into $V\subset W$.

\hfill

{\bf Step 1:} Let $x, y \in \Perspace$ be distinct
points, and $V_x, V_y\subset H^2(M, \R)$ the associated
2-planes. Then $x$ and $y$ belong to the same
hyperk\"ahler line $S$ if and only if $V_x\cap V_y$ is non-zero,
and the space $\langle V_x, V_y\rangle$
generated by $V_x, V_y$ is positive. This is 
an immediate consequence of Step 0.

\hfill

{\bf Step 2:} Let $x\in \Perspace$ be a 
point, and $V_x\subset H^2(M, \R)$ the corresponding 2-plane.
A vector $\omega\in V_x^\bot$  in the positive cone of
$V_x$ defines 
a 3-dimensional plane $\langle \omega, V_x\rangle$.
This gives a hyperk\"ahler line $C_\omega\subset \Teich$
passing through $x$, whenever $q(\omega, \omega)>0$.  Clearly,
for generic $\omega \in V_x^\bot$, all rational points of $\omega^\bot$
lie in $(H^{2,0}\oplus H^{0,2})\cap H^2(M, \Q)$. 
Therefore, the orthogonal complement to 
$H^{2,0}\oplus H^{0,2}\oplus \R\omega$ has no rational points
(see also \ref{_GHK_lines_and_Pic_Claim_}).

\hfill

{\bf Step 3:} Let $W_1$ and $W_2$ be 3-dimensional
positive subspaces in the space $H^2(M, \R)$,
containing $a\in  H^2(M, \R)$. Assume that $a^\bot\cap H^2(M, \Q)=0$.
By \ref{_GHK_lines_and_Pic_Claim_}
this implies, in particular, that the subspaces 
$W_i$ correspond to  GHK lines $S_{W_1}$, $S_{W_2}$. 
Then there exists a GHK line intersecting $S_{W_1}$ 
and $S_{W_2}$. Indeed,
from \ref{_vector_posi_def_Lemma_} it follows that
there exists a positive 2-dimensional plane 
$V:=\langle a, z \rangle\subset H^2(M, \R)$ generated by
$a, z$, with $z\in W_1$. Applying \ref{_vector_posi_def_Lemma_},
again, we find a positive 3-dimensional plane 
$W:=\langle a, z, z'\rangle\subset H^2(M, \R)$,  
with $z' \in W_2$. By \ref{_GHK_lines_and_Pic_Claim_}, 
$W$ corresponds to a GHK line $S_W$. Now, Step 1 immediately
implies that $S_W$ intersects $S_{W_1}$ and $S_{W_2}$.

\hfill

{\bf Step 4:} Let  $x, t \in \Perspace$.
Using Step 2, we find
GHK lines passing through $x$ and $t$. Denote
by $W_x$, $W_t\subset H^2(M, \R)$ the corresponding 3-planes,
and let $a\in W_x$ be a vector which 
satisfies  $a^\bot \cap H^2(M, \Q)=0$ 
(such $a$ exists by \ref{_GHK_lines_and_Pic_Claim_}.)
Using \ref{_vector_posi_def_Lemma_}, 
we choose a non-zero $b\in W_t$, in such a way that
$q|_{\langle a,b\rangle}$ is positive definite. Now, let  
$W$ be a positive 3-plane in $H^2(M, \R)$ containing
$a$ and $b$. By Step 3, there exist a GHK line
intersecting $S_{W_x}$ and $S_W$, and another
GHK line intersecting $S_{W_t}$ and $S_W$.
We proved  \ref{_4_GHK_lines_Proposition_}.
\endproof


\section{GHK lines and exceptional sets}
\label{_HK_lines_Section_}


\subsection{Lifting the GHK lines to the Teichm\"uller space}

The following proposition ensures that
GHK lines are in some sense 
``lift\-able'' to the Teichm\"uller space.
This is a key idea used to prove that
the period map is a covering.

\hfill

\proposition\label{_lifting_GHK_lines_Proposition_}
Let $I\in \Teich$ be a point in the Teichm\"uller space
of a hyperk\"ahler manifold, $\NS(M,I)=0$, and $S\subset
\Perspace$ a hyperk\"ahler line passing through 
$\Per(I)$.\footnote{Such a hyperk\"ahler line is
  necessarily generic, by \ref{_GHK_lines_and_Pic_Claim_}.} 
Then there exists a holomorphic curve $S_I\subset \Teich$
passing through $I$ and satisfying $\Per(S_I)=S$.

\hfill

{\bf Proof:} Denote by $W\subset H^2(M, \R)$ 
the 3-dimensional space used to define $S$.
Let $\Omega$ be the holomorphic symplectic
form of $(M,I)$, and 
$V := \langle \Re\Omega,\Im\Omega\rangle\subset H^2(M, \R)$
the corresponding 2-dimensional space.
Then $V\subset W$, and the 1-dimensional
orthogonal complement $V^\bot\cap W$ intersects both
components of the cone $\{x\in H^{1,1}_I(M, \R)\ \ |\ \ q(x,x)>0\}$.
One of these components coincides with the K\"ahler cone
(\ref{_Kahler_cone_tri_NS_Corollary_}). 
Choose a K\"ahler form $\omega \in V^\bot_W$,
normalize it in such a way that 
\[ 
  q(\Re\Omega, \Re\Omega) = q(\Im\Omega,\Im\Omega)= q(\omega,\omega),
\]
and let $(M,I,J,K)$ be the hyperk\"ahler structure
associated with $\omega$ as in \ref{_Calabi-Yau_Theorem_}.
Denote by $S_I$ the line of complex structures
associated with this hyperk\"ahler structure.
As shown above (\ref{_sphere_hk_Remark_}), the period map
$\Per$ maps $S_I$ isomorphically
to $S$. 
\endproof

\hfill

Abusing the language, we call 
a $\C P^1$ of induced complex structures associated
with a hyperk\"ahler structure ``a hyperk\"ahler line'' as well.
These ``hyperk\"ahler lines'' lie in the Teichm\"uller space,
and the hyperk\"ahler lines defined previously lie in
the period space. Then \ref{_lifting_GHK_lines_Proposition_}
can be restated saying that a GHK line passing through
a point $l\in \Perspace$, satisfying $\NS(M,l)=0$,
can be always lifted to a hyperk\"ahler line
$S\subset \Teich$ for each $I\in \Teich$
such that $\Per(I)=l$.

\hfill

\definition
Let $\Perspace$ be a period space for a hyperk\"ahler manifold $M$,
and $\psi:\; D \arrow \Perspace$ an etale map from a Hausdorff 
manifold $D$. Given a hyperk\"ahler line $S\subset \Perspace$,
denote by $S_{\Pic >0}$ the set of all $I\in S$
satisfying $\rk \Pic(M,I)>0$. We say that
$\psi$ is {\bf compatible with generic hyperk\"ahler lines}
if for any GHK line $S\subset \Perspace$,
the space $X:= \psi^{-1}(S)$ is a union of several 
disjoint copies of $S$, which are closed and open in $X$, and another
subset $Y\subset X$, which satisfies
$\psi(Y)\subset S_{\Pic >0}$.

\hfill

\proposition\label{_period_compa_GHK_l_Proposition_}
Let $M$ be a hyperk\"ahler manifold, and 
\[ \Per:\; \Teich_b \arrow \Perspace\]
its period map. Then $\Per$ is compatible with
generic hyperk\"ahler lines.

\hfill

{\bf Proof:} Let $S\subset \Perspace$ be a GHK line, 
$l\in S$ a point with $\NS(M,l)=0$, and $I\in \Teich$
a point in the fiber $\Per^{-1}(l)$. By \ref{_lifting_GHK_lines_Proposition_},
$S$ can be lifted to a hyperk\"ahler line $S_I\subset \Teich$
passing through $I$. Since $\Per$ is etale,
the restriction $\Per:\; S_I\arrow S$
is a diffeomorphism. By \ref{_preima_isola_compo_Claim_} below,
$S_I$ is a connected component of 
$\Per^{-1}(S)$.
\endproof

\hfill

The following claim is completely trivial.

\hfill

\claim\label{_preima_isola_compo_Claim_}
Let $X \stackrel \psi \arrow Y$ be a local homeomorphism
of Hausdorff spaces, 
$S\subset Y$ a compact subset, and $S_1\subset X$ a subset
of $\psi^{-1}(S)$, with 
$\psi\restrict{S_1}:\; S_1 \arrow S$
a homeomorphism. Then $S_1$ is closed and open
in $\psi^{-1}(S)$

\hfill

{\bf Proof:} The set $S_1$ is closed because
it is homeomorphic to $S$ which is compact, and $X$ is Hausdorff.
Suppose that 
$S_1$ is not open in $\psi^{-1}(S)$; then, there exists a sequence
of points $\{x_i\}\subset \psi^{-1}(S)\backslash S_1$ converging
to $x\in S_1$. Choose a neighbourhood $U\ni x$
such that $\psi\restrict U$ is a homeomorphism.
Replacing $\{x_i\}$ by a subsequence, 
we may assume that $\{x_i\}\subset U$.
Then $\psi \restrict{S_1 \cap U}$ is a homeomorphism
onto its image $S_U$, which is a neighbourhood of 
$\psi(x)$ in $S$. Replacing $\{x_i\}$ by a subsequence
again, we may assume that all $\psi(x_i)$ lie in
$S_U$. Since $\psi \restrict{U}$ is bijective onto
its image, this map induces
a bijection from  $S_1 \cap U$ to $S_U$.
Therefore, $\{x_i\}\subset S_1\cap U$.
We obtained a contradiction, proving that
$S_1$ is open in $\psi^{-1}(S)$.
 \endproof

\subsection{Exceptional sets of etale maps}

In \cite{_Browder_}, F. Browder
has discovered several criteria which can be used
to prove that a given etale map is a covering.
Unfortunately, in our case neither of his
theorems can be applied, and we are forced to
devise a new criterion, which is then
applied to the period map.

\hfill

\definition
Let $X\stackrel \psi\arrow Y$ be a local homeomorphism 
of Hausdorff topological spaces, {\em e.g.} an etale map.
Consider a connected, simply connected subset  $R\subset Y$,
and let $\{R_\alpha\}$ be the set of connected components
of $\psi^{-1}(R)$.  {\bf An exceptional set} of $(\psi, R)$
is $R \backslash \psi(R_\alpha)$. 

\hfill

\remark
The following topological criterion is the main 
technical engine of this section. Its proof is
complicated, but completely abstract, and
we hope that this result might have independent
uses outside of hyperk\"ahler geometry.
We include an alternative proof of this proposition 
in the Appendix by Eyal Markman (Section \ref{_Appe_Section_}).

\hfill

\proposition\label{_exce_empty_covering_Proposition_}
Let $X\stackrel \psi\arrow Y$ be a local diffeomorphism
of Hausdorff manifolds.
Assume that for any open subset $U\subset Y$,
the closure $\bar U \subset Y$ has empty exceptional sets,
provided that $U$ has smooth boundary. Then 
$\psi$ is a covering.

\hfill

{\bf Proof:} \ref{_exce_empty_covering_Proposition_}
is local in $Y$, hence it will suffice to prove it when $Y$ is
diffeomorphic to $\R^n$. Choose a flat Riemannian
metric on $Y \cong \R^n$. Lifting the corresponding
Riemannian metric to $X$, we can consider $X$
as a Riemannian manifold, also flat. The Riemannian
structure defines a metric on $Y$ and $X$ as
usual. For a point $x$ in a metric space $M$, 
{\bf a closed $\epsilon$-ball 
with center in $x$} 
is the set 
\[
\bar B_\epsilon(x) := \{ m \in M \ \ | \ \ d(x, m) \leq \epsilon\}.
\]
Taking strict inequality, we obtain {\bf an open ball},
\[
B_\epsilon(x) := \{ m \in M \ \ | \ \ d(x, m) < \epsilon\}.
\]

Clearly, $\bar B_\epsilon(x)$ is closed,
$B_\epsilon(x)$ is open, and $\bar B_\epsilon(x)$ 
is the closure of $B_\epsilon(x)$, and its completion,
in the sense of metric geometry.

For any
$x\in X$, $y=\psi(x)$, let $D_x \subset \R^{>0}$ be 
the set of all $\epsilon \in \R^{>0}$
such that the corresponding $\epsilon$-ball 
$\bar B_\epsilon(x)$ is mapped
to $\bar B_\epsilon(y)$ bijectively.
Clearly, $D_x$ is an initial interval of $\R^{>0}$.
We are going to show that $D_x$ is open and
closed in $\R^{>0}$.

\hfill

{\bf Step 1:} The interval $D_x$ 
is open, for any etale map
$X\stackrel \psi\arrow \R^n$. Indeed,
for any $\epsilon \in D_x$,
the corresponding $\epsilon$-ball
$\bar B_\epsilon(x)$ is compact, because
it is isometric to $\bar B_\epsilon(y)$.
Every point of $\bar B_\epsilon(x)$
has a neighbourhood which is
isometrically mapped to 
its image in $Y$. Take a 
covering $\{B_\epsilon(x), U_1, U_2, ...\}$ of $\bar B_\epsilon(x)$ 
where  $U_i$ are open balls with this property, centered in a 
point on the boundary 
of $\bar B_\epsilon(x)$. Since $\bar B_\epsilon(x)$
is compact, $\{B_\epsilon(x), U_i\}$ has a finite subcovering $U_1, ..., U_n$.
By construction, for each point 
$z\in W:= B_\epsilon(x)\cup 
\bigcup_i U_i$, the set $W$ contains
a straight line (geodesic) from $x$ to $z$. Indeed, 
$W$ is a union of an open ball $B_\epsilon(x)$
and several open balls centered on its boundary, 
and all these balls are isometric to open
balls in $\R^n$.  Since $\psi$ maps
straight lines to straight lines, it maps
$B_{\epsilon'}(x)$ surjectively to $B_{\epsilon'}(y)$.
To show that this map is also injective, 
consider two points $a_1, a_2 \in B_{\epsilon'}(x)$,
mapped to $b\in B_{\epsilon'}(y)$, and let
$[x, a_1]$ and $[x, a_2]$ be the corresponding
intervals of a straight line.
Since $\psi(a_1)=\psi(a_2)=b$, one has $\psi([x, a_1]) = \psi([x, a_2])$, and
these intervals have the same length. Also, 
$[x, a_1]\cap B_\epsilon(x)=[x, a_2]\cap B_\epsilon(x)$,
because $\psi\restrict{B_\epsilon(x)}$ is injective.
Therefore, the intervals  $[x, a_1]$ and $[x, a_2]$
coincide, and $a_1=a_2$.

\hfill

{\bf Step 2:}
Let $\psi:\; X \arrow \R^n$ be an etale map,  
$y=\psi (x)$, and suppose that $\psi:\;B_s(x)\arrow  B_s(y)$
is bijective, for some $s>0$. Then $\phi:\; B_s(x) \arrow B_s(y)$
is an isometry, with respect to the metric on $B_s$ induced 
from the ambient manifold. Indeed, $\psi$ is etale, hence
any piecewise geodesic path in $X$ is projected to such
one in $\R^n$. Therefore, $\psi$ does not increase distance:
$d(a, b) \geq d(\psi(a), \psi(b))$. The open ball
$B_s(y)$ is geodesically convex, hence 
for any $y_1, y_2 \in B_s(y)$, the geodesic interval
$[y_1, y_2]$ can be lifted to a geodesic in $B_s(x)$.
This implies an inverse inequality: $d(a, b) \leq d(\psi(a), \psi(b))$.
We proved that $\phi:\; B_s(x) \arrow B_s(y)$ is an isometry.
This implies that the map 
 $\psi:\;\bar B_s(x)\arrow \bar  B_s(y)$ of their metric completions
is also an isometry. In particular, this map is injective.

\hfill

{\bf Step 3.} In the assumptions of 
Step 2, we prove that $\bar B_s(x)$
is a connected component of $\psi^{-1}(\bar B_s(y))$.
Notice that $\bar B_s(x)$ is a closure of
$B_s(x)$, which is homeomorphic to a ball in $\R^n$,
hence $\bar B_s(x)$  is connected. To prove
that it is a connected component, we need only
to show that it is open in $\psi^{-1}(\bar B_s(y))$.

The corresponding
map of open balls $\psi:\;B_s(x)\arrow  B_s(y)$
is by definition bijective. The closed ball
$\bar B_s(x)$ is closed in $\psi^{-1}(\bar B_s(y))$. For any 
$z\in\6\bar B_s(x)$ on the boundary of
$\bar B_s(x)$, an open ball $S$ centered
in $z$ is split by the boundary 
\[ \6\bar B_s(x)=\{x'\in X\ \ | \ \ d(x,x')=s\}
\]
onto two open components, 
$S_o:=\{x'\in X\ \ | \ \ d(x,x')>s\}$ and 
$S_1:=\{x'\in X\ \ | \ \ d(x,x')<s\}$, with 
$S_1$ mapping to $B_s(y)$, $\6\bar B_s(x)$
mapping to its boundary, and $S_o$
to $Y\backslash \bar B_s(y)$.\footnote {Here we use
the fact that $\psi\restrict S$ is a bijection,
for $S$ sufficiently small, hence the image of 
$S$ cannot wrap on itself.} This implies that
\[ 
 \psi^{-1}(\bar B_s(y)\cap \psi(S))=S\cap \bar B_s(x).
\]
Therefore, 
$\bar B_s(x)$ is open in $\psi^{-1}(\bar B_s(y))$.

\hfill

\centerline{\scalebox{0.7}{\includegraphics{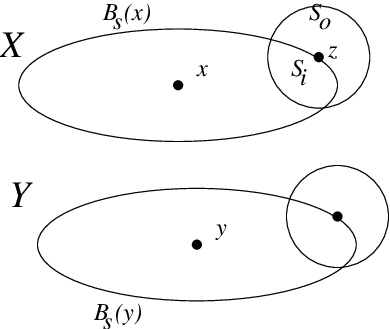}}}

\hfill

{\bf Step 4.} Now we can show that $D_x$ is closed.
This argument uses the triviality of exceptional sets
(the first time in this proof, the rest follows just 
from the  etaleness of $\psi$).
Let $s:= \sup D_x$, and $\bar B_s(x)$
the corresponding closed ball. 
We prove that $\psi:\;\bar B_s(x)\arrow \bar  B_s(y)$
is a homeomorphism.

From Step 3, it follows that $\bar B_s(x)$ is a 
connected component of the preimage $\psi^{-1}(\bar B_s(y))$.
Since the exceptional sets of $\bar B_s(y)$
are all empty, the map $\psi:\;\bar B_s(x)\arrow \bar  B_s(y)$
is surjective. It is injective as follows from Step 2.

We proved that $D_x$ is open and
closed, hence $D_x= \R^{>0}$, and $\psi$ maps
any connected component of $X$ bijectively to $Y$.
\endproof

\hfill

\remark
An exceptional set of $(\psi, U)$ is always closed in $U$.

\hfill 

\lemma\label{_exce_covered_curves_Lemma_}
Let $M$ be a Hausdorff manifold, 
$M \stackrel \psi \arrow \Perspace$ a local
diffeomorphism, compatible with GHK lines, 
$U\subset \Perspace$ an
open, simply connected subset, 
$U_\alpha$ a component of $\psi^{-1}(U)$,
and $K_\alpha$ the corresponding exceptional set. 
Consider a GHK line $C\subset \Perspace$, and 
let $C_1$ be a connected component of $C\cap U$.
Then $C_1\subset K_\alpha$, or 
$C_1 \cap K_\alpha=\emptyset$.

\hfill

{\bf Proof:} Suppose that 
$D:= C_1 \cap (U\backslash K_\alpha)$
is non-empty. Since $K_\alpha$ is closed in $U$,
$D$ is open in $C_1$. Then $D$ contains points $l\in D$
with $\NS(M,l)=\emptyset$ (\ref{_gene_on_hk_line_Remark_}). 
The set $\psi^{-1}(l)$ is non-empty, because $l\notin K_\alpha$.
Since $\psi$ is compatible with GHK lines,
for any $I \in \psi^{-1}(l)$, there is a curve 
$C_I\subset M$ passing through $I$ and projecting
bijectively to $C$. Clearly, the connected 
component of $C_I \cap \psi^{-1}(U)\ni I$ is bijectively mapped to
$C_1$, hence $C_1 \cap K_\alpha=\emptyset$.
\endproof

\hfill

\remark\label{_exce_covered_curves_boundary_Remark_}
A version of \ref{_exce_covered_curves_Lemma_} 
is also true if $\bar U$ is a closed set, obtained
as a closure of an open subset $U \subset \Perspace$,
and $C_1$ a connected component of $\bar U \cap C$,
for a GHK curve $C$. If $C_1$ contains interior
points, the same argument as above can be used
to show that $C_1\subset K_\alpha$, or 
$C_1 \cap K_\alpha=\emptyset$.

\subsection{Subsets covered by GHK lines}

Let $U\subset \Perspace$ be an open subset, 
or a closure of an open subset with smooth boundary, 
and $K\subset U$ a subset of $U$. 
Given a GHK line $C\subset \Perspace$, denote by
$C_U$ a connected component of $C\cap U$.
This component is non-unique
for some $C$ and $U$.
Denote by $\Omega_U(K)$ the union of all segments
$C_U\subset U$ intersecting $K$,
for all GHK lines $C\subset \Perspace$.
In other words, $\Omega_U(K)$ is the
set of all points connected to $K$ by a 
connected segment
of $C \cap U$, with $C\subset \Perspace$ a GHK line.
Let $\Omega^*_U(X)$ be the union of
$\Omega_U(X), \Omega_U(\Omega_U(X)), \Omega_U(\Omega_U(\Omega_U(X))), ...$.

\hfill

\proposition\label{_Omega_*_open_Proposition_}
Let $U\subset \Perspace$ be an open subset,
and $x\in U$ a point. Then $\Omega_U^*(x)$
is open in $U$.

\hfill

{\bf Proof:} 
We give two alternative proofs of this proposition.
The first one is based on \cite{_Beau_et_al:K3_} and
\cite[Theorem 3.1, Theorem 5.1]{_Verbitsky:coho_announce_},
where essentially the same argument was used to show that the category of
polystable holomorphic vector bundles on a generic hyperk\"ahler
manifold $M$ is independent from the choice of $M$
in its deformation class. A beautiful presentation of
this argument is found in 
\cite[Proposition 3.7]{_Huybrechts:Bourbaki_}.

The second proof uses the concept of subtwistor metric $d_{tw}$,
introduced in the Appendix, Section \ref{_Appe_subtwistor_}.
This argument, based on the theory of intrinsic metrics
and the Gleason-Palais's elaboration of the
Gleason-Montgomery-Zippin solution of Hilbert's
fifth problem, is much more explicit.

{\bf Proof of \ref{_Omega_*_open_Proposition_}: the
first approach.} From 
\ref{_4_GHK_lines_Proposition_}, it follows that 
$\Omega_U^5(x)$ is open. Indeed,
any $y$ can be connected to $x$ by a sequence ${\goth S}(x,y)$
of 5 GHK lines. From its proof it is apparent that 
one can chose this sequence in such a
way that it depends continuously on $x,y$. 
Let $y\in \Omega_U^5(x)$, and let 
$\gamma:\; [0,1]\arrow U$ be a path connecting
$x$ to $y$ and sitting in the union of the GHK lines 
$S_i\in {\goth S}(x,y)$. Given a sequence $\{y_i\}$ converging to $y$,
and a sequence ${\goth S}(x,y_i)$ converging to ${\goth S}(x,y)$,
we can chose a sequence of paths $\gamma_i:\; [0,1]\arrow \Perspace$
with the following properties.

\begin{itemize}
\item The paths $\gamma_i$ sit in ${\goth S}(x,y_i)$.
\item The sequence $\{\gamma_i\}$ converges to $\gamma$
in the compact-open topology.
\item Each $\gamma_i$ connects $x$ to $y_i$.
\end{itemize}

Since $[0,1]$ is compact, $U$ is open, and $\gamma([0,1])\subset U$,
for any sufficiently big $i$, one has $\gamma_i([0,1])\subset U$.
This implies that $y_i\in \Omega^5_U(x)$.

{\bf Proof of \ref{_Omega_*_open_Proposition_}: the
second approach.} 
It suffices to prove \ref{_Omega_*_open_Proposition_}
when the closure of $U$ is compact. Indeed, each
point $z\in \Omega^*_U(x)$ has a neighbourhood $V\subset U$
with its closure in $U$ compact. If $\Omega^*_V(x)$
is open in $V$, then each point $z\in \Omega^*_U(x)$
has an open neighbourhood $V_1\subset V\subset U$ 
contained in $\Omega^*_U(x)$.

By \ref{_d_tw_metric_Theorem_} (see the Appendix,  
Section \ref{_Appe_subtwistor_}), the metric 
$d_{tw}$ induces the usual topology on $\Perspace$.
For any $x\in U$, the distance $d_{tw}(x, \6 U)$ to the
boundary of $U$ is positive,
because $\6 U$ is compact.
Then, for any $r<d_{tw}(x, \6 U)$, 
the open ball $B_r(x, d_{tw})$ is contained in $\Omega_U^*(x)$.
Indeed, let $y\in B_r(x, d_{tw})$. Then $y=s_{n+1}$
is connected to $x=s_0$ by a sequence of
GHK lines $S_1, ..., S_n$, such that
$\sum_{i=0}^n d_g(s_i, s_{i+1})< r$.
Consider the corresponding piecewise
geodesic path $\gamma\subset \bigcup S_i$ of length $<r$.
Since $d_{tw}(x, \6 U)>r$, the whole of $\gamma$
belongs to $U$. Therefore, $y$ is connected to 
$x$ by a union of connected segments of GHK lines
which lie in $U$. 
\endproof

\hfill

To apply \ref{_exce_empty_covering_Proposition_} 
to the period map using the exceptional sets,
we also need closed subsets with smooth boundary.
In this situation the following lemma  can be used.

\hfill

\lemma\label{_Omega_with_bou_Lemma_}
Let $K\subset \Perspace$ be a compact closure of an open subset
with smooth boundary, and $x\in K$ a point.  Then $\Omega_K(x)$
contains an interior point of $K$.

\hfill

{\bf Proof:} 
Let $V_x$ be the 2-plane
in $H^2(M, \R)$ corresponding to $x$ via the
identification $Gr(2)=\Perspace$. Then
the tangent space $T_x\Per$ is identified
with $\Hom(V_x, V_x^\bot)$, where
$V_x^\bot$ is the orthogonal complement.
For a hyperk\"ahler line $C$ associated with
a 3-dimensional space $W$, the corresponding
2-dimensional space $T_x C\subset T_x\Perspace$
is the space $\Hom(V_x, (V_x^\bot\cap W))$.
Since $V_x^\bot=H^{1,1}_x(M)$ and 
$W$ can be chosen by adding 
to $V_x$ any K\"ahler class in $H^{1,1}(M)$,
the set of all tangent vectors $T_x C\subset T_x\Perspace$
is open in the space
\[
P:=\{ l \in \Hom(V_x, H^{1,1}_x(M)) \ \ | \ \ \rk l =1 \}
\]
The condition $\rk l =1$ is quadratic, and
it is easy to check that an open subset $U_P\subset P$ cannot be 
contained inside a linear subspace of 
positive codimension. In particular,
$U_P$ cannot lie in the tangent space to the boundary of $K$,
\begin{equation}\label{_U_P_not_tangent_Equation_}
 U_P\not\subset T_x\6 K\subset\Hom(V_x, H^{1,1}(M)).
\end{equation}
Take for $U_P$ the set of all vectors tangent to
GHK lines passing through $x$. Then 
\eqref{_U_P_not_tangent_Equation_} implies that
for a generic
GHK line $C$ passing through $x$, $C$
intersects with the interior points of $K$.
 \endproof

\hfill

\corollary\label{_Omega_all_Corollary_}
Let $K\subset \Perspace$ be 
a closure of an open, connected subset $U\subset \Perspace$
with smooth boundary,
and $\Omega_K$ the operation on subsets of $K$ defined above.
Then $\Omega_K^*(x)=K$, for any point $x\in K$.

\hfill

{\bf Proof:} 
Clearly, $\Omega_U^*(x)$ is the set of all points in $U$
which can be connected to $x$ within $U$ by a finite sequence
of connected segments of GHK lines. By 
\ref{_Omega_*_open_Proposition_}, 
$\Omega_U^*(x)$ is open in $U$.   
If $y\notin \Omega_U^*(x)$, then $\Omega_U^*(y)$
does not intersect $\Omega_U^*(x)$. Then
$U$ is represented as a disconnected union
of open sets $\Omega_U^*(x_i)$, for some $\{x_i\}\subset U$.
This is impossible, because $U$ is connected.
We proved that $\Omega^*_U(x)=U$.
Then $\Omega_K^*(x)=K$, because every point
on a boundary of $K$ is connected to some
point of $U$ by a connected segment of
a GHK line (\ref{_Omega_with_bou_Lemma_}).
\endproof

\hfill

The main result of this section is the following theorem

\hfill

\theorem\label{_curves_then_covering_Theorem_}
Let $M$ be a Hausdorff manifold and
$M\stackrel\psi\arrow \Perspace$ a local
diffeomorphism compatible with GHK lines. Then 
$\psi$ is a covering.

\hfill

\remark
It is well known that $\Perspace$ is simply connected
(\ref{_perspace_simply_conne_Claim_}).
Then \ref{_curves_then_covering_Theorem_}
implies that $\psi$ is a diffeomorphism.

\hfill

{\bf Proof of \ref{_curves_then_covering_Theorem_}:}
To prove that $\psi$ is a covering,
it suffices to show that all its exceptional sets 
of $(\psi, K)$, are empty provided that
$K$ is a closure of a simply connected open
subset $U\subset \Perspace$ which has a smooth boundary
(\ref{_exce_empty_covering_Proposition_}). 
Let $K_\alpha$ be an exceptional set, associated
with a closure $K\subset \Perspace$ of an
open subset $U\subset \Perspace$ with smooth boundary.
From \ref{_exce_covered_curves_Lemma_} and
\ref{_exce_covered_curves_boundary_Remark_}
it follows that $\Omega_K(K_\alpha)= K_\alpha$,
where $\Omega_K(Z)$ is a union of all
connected segments of $C\cap K$ intersecting
$Z$, for all GHK lines $C\subset \Perspace$. 
Then $\Omega_K^*(K_\alpha)=K_\alpha$,
where $\Omega_K^*(Z)$ is a union of 
all iterations $\Omega_K^i(Z)$. However, for 
any non-empty $Z\subset K$, one has
$\Omega_K^*(Z)=K$ by \ref{_Omega_all_Corollary_}.
Therefore, any exceptional set $K_\alpha$
of $(\psi,K)$ for $K$ as above is empty, and 
\ref{_curves_then_covering_Theorem_} follows.
\endproof

\hfill

\remark\label{_Per_cove_Proposition_proof_Remark_}
\ref{_Per_covering_Proposition_}
is implied by \ref{_period_compa_GHK_l_Proposition_} 
and \ref{_curves_then_covering_Theorem_} below.
Indeed, by \ref{_period_compa_GHK_l_Proposition_}, 
$\Per$ is compatible with the generic hyperk\"ahler lines 
(\ref{_GHK_line_Definition_}),
and by \ref{_curves_then_covering_Theorem_}, 
any such map is necessarily a covering.

\hfill


\section{Monodromy group for $K3^{[n]}$.}


 When $M=K3^{[n]}$ is a Hilbert scheme of points
on a K3 surface, fundamental results about its moduli
were obtained by E. Markman (\cite{_Markman:mono_}, 
\cite{_Markman:constra_}), using the Fourier-Mukai 
action on the derived category of coherent sheaves.
In this section we relate these results with our
computation of $\Teich_b$ to obtain a global
Torelli theorem for $M=K3^{[p^\alpha+1]}$, $p$ prime.

\subsection{Monodromy group for hyperk\"ahler manifolds}
\label{_mono_hype_Subsection_}

\definition
The {\bf monodromy group} $\Mon(M)$ 
of a hyperk\"ahler manifold
$M$ is a subgroup of $GL(H^*(M, \Z))$ generated by 
the monodromy of the Gauss-Manin local systems, for all
holomorphic deformations of $M$ over a 
connected complex analytic base.

\hfill

Using the global Torelli theorem
(\ref{_global_Torelli_Theorem_}), 
the monodromy group can be related to the
mapping class group, as follows.

\hfill

\theorem\label{_mono_via_mappi_Theorem_}
Let $(M,I)$ be a hyperk\"ahler manifold, 
and $\Teich^I$ the corresponding connected
component of a Teichm\"uller space. 
Denote by $\Gamma^I$ the subgroup
of the mapping class group preserving 
the component $\Teich^I$, and let $\Mon(M)$ be the
monodromy group of $(M,I)$ defined above.
Then  $\Mon(M)$ coincides
with the image of $\Gamma^I$ in $GL(H^*(M, \Z))$.

\hfill

{\bf Proof:}\footnote{I am grateful to the referee for numerous
suggestions which lead to many improvements in this proof.}
From its construction, it is obvious that 
the Teichm\"uller space is the coarse moduli space for the
following functor  from analytic spaces to sets. 
This functor associates to a complex analytic space $B$
the set of isomorphism classes of pairs $({\cal X}\arrow B,\Psi)$, 
where ${\cal X}\arrow B$ is a complex analytic deformation 
of $M$ over $B$, and $\Psi$ a smooth trivialization
of the family ${\cal X}\arrow B$, defined up to isotopy
on the fibers. 

Now, consider an element $\gamma\in\Mon(M)$
in the image of the monodromy of a holomorphic family
${\cal Z}\arrow B$. Then there exists
a point $b\in B$ and a loop $\gamma_0:[0,1]\arrow B$ such that
the corresponding Gauss-Manin monodromy induces
$\gamma$ on $H^*(M,\Z)$. 

For a certain covering
$\tilde B \arrow B$, the corresponding family
$\tilde {\cal Z}\arrow\tilde  B$ admits a
smooth trivialization. This gives a map
$\tilde  B\stackrel \psi \arrow \Teich^I$ to the corresponding
coarse moduli space, such that the pullback $\tilde {\cal Z} \arrow \tilde B$
of the family ${\cal Z}\stackrel \pi \arrow B$ admits a smooth trivialization 
$\tilde {\cal Z}=\tilde B\times M$. 

Let $\tilde \gamma_0$ be 
a lifting of $\gamma_0$ to $\tilde B$, and $x,y\in \tilde B$
the ends of this path, with $\tilde {\cal Z}_x=(M,I)=\tilde {\cal Z}_y$
denoting the fibers of $\pi$ over $x,y$. The trivialization
of $\pi$ over $\tilde B$ induces a diffeomorphism
$(M,I)=\tilde {\cal Z}_x\arrow \tilde {\cal Z}_y=(M,I)$
acting on $H^*(M,\Z)$ as $\gamma$. 

In a neighbourhood of the corresponding path 
$\psi(\gamma_0)\subset  \Teich^I$ the universal family 
of deformations of $M$ is well defined. Replacing $B$ by
a smaller neighbourhood of $\gamma_0$ if necessary,
we may assume that the family 
$\tilde {\cal Z} \arrow \tilde B$ is a pullback
of the universal family on a neighbourhood $U$ of 
$\psi(\gamma_0)$. The monodromy group of this family
by definition belongs to $\Gamma_I$, hence the image of
$\gamma$ in $GL(H^*(M,\Z)$ lies in $i(\Gamma_I)$.

Conversely, for each $\gamma \in \Gamma^I$,
consider the action of $\Gamma^I$ on $\Teich^I$,
and let $x, y:= \gamma(x)$ be a pair of points on $\Teich^I$,
connected by a smooth path $\gamma_0$. Denote by
$U$ a neighbourhood of $\gamma_0$, diffeomorphic 
to an open ball. Consider a non-normal quotient $U_1$,
obtained from $U$ by identifying $x$ and $y$.
Since $U$ is diffeomorphic to an open ball,
there exists a universal fibration ${\tilde X}\arrow
U$. 
Gluing two fibers of this universal fibration, we obtain
a holomorphic fibration ${\tilde X}_1\arrow U_1$; its monodromy acts
on $H^2(M,\Z)$ as $\gamma$, by construction. 
 \endproof

\hfill

This result allows one to answer the question
asked in \cite{_Markman:constra_} (Conjecture 1.3). 

\hfill

\corollary\label{_Markman_question_Corollary_} 
Let $\gamma\in \Mon$ be an element of the
monodromy group acting trivially on 
the projectivization ${\Bbb P}H^2(M,\C)$.
Suppose that a general deformation of $M$ has no 
automorphisms. Then $\gamma$ is trivial.

\hfill

{\bf Proof:} Using \cite[Remark 13]{_Catanese:moduli_}, 
we may assume that there exists a universal fibration 
${\cal Z}\stackrel \pi \arrow \Teich_I$. This gives a local system 
$R\pi_*\Z$ over $\Teich^I$. Let $\gamma$ be 
an element of the mapping class group acting trivially on
$H^2(M)$ preserving a connected component $\Teich^I\subset \Teich$.
By \ref{_stabili_in_H^2_maping_class_Theorem_}, $\gamma$
acts trivially on $\Teich^I$. However, \ref{_mono_via_mappi_Theorem_}
implies that an action of any $\gamma\in \Mon$ is induced by
the parallel transport in the local system $R\pi_*\Z$
from $x$ to $\gamma(x)$ along a path $\gamma_0$
connecting $x$ to $\gamma(x)$ in $\Teich^I$. Since the action
of $\gamma$ on $\Teich^I$ is trivial, we may choose $\gamma_0$ 
to be trivial.
\endproof

\hfill

\remark
In the above corollary, a stronger result is actually
proven. Instead of defining the monodromy group as above,
we could define $\tilde\Mon$ as the image
of $\pi_1({\cal M})$ in the mapping class
group of $M$. Then \ref{_Markman_question_Corollary_}
implies that the natural map of
$\tilde \Mon$ to $PGL(H^2(M, \C))$
is injective.

\hfill

\remark
The kernel of the natural projection 
$\Gamma_I\arrow PGL(H^2(M, \C)$ is identified with the
group of holomorphic automorphisms of a generic
deformation of a hyperk\"ahler manifold $M$.
When $M$ is a deformation of 
a Hilbert scheme of a K3, this group is trivial,
which can be easily seen e.g. from the results of
\cite{_V:K3^[n]_}. When $M$ is a generalized
Kummer variety, it is known to be
non-trivial (\cite{_KV:Kummer_Dynkin_}).

\subsection{The Hodge-theoretic
Torelli theorem for $K3^{[n]}$}
\label{_Hodge_Torelli_Subsection_}
\newcommand{\Ref}{\operatorname{Ref}}

\definition
Let $V$ be a vector space, $g$ a non-degenerate quadratic
form, and $v\in V$ a vector which satisfies
$g(v, v) = \pm 2$. Consider the 
pseudo-reflection map $\rho_v:\; V \arrow V$,
\[
\rho_v(x) := \frac{-2}{g(v,v)} x + g(x, v) v.
\]
Clearly, $\rho_v$ is a reflection when $g(v, v) = 2$,
and $-\rho_v$ is a reflection when $g(v,v)=-2$.
Given an integer lattice in $V$, consider the group
$\Ref(V)$ generated by $\rho_v$ for all integer vectors
$v$ with $g(v, v) = \pm 2$. We call $\Ref$ {\bf a reflection
group}.

\hfill

The following fundamental theorem was proven by
E. Markman in \cite{_Markman:constra_}.

\hfill

\theorem
(\cite[Theorem 1.2]{_Markman:constra_})
Let $M=K3^{[n]}$ be a Hilbert scheme of points on a 
K3 surface,  and $\Mon^2$ be the image of the monodromy
group in $GL(H^2(M,\Z))$. Then
$\Mon^2= \Ref(H^2(M, \Z),q)$.
\endproof

\hfill

Comparing this with \ref{_mono_via_mappi_Theorem_} 
and using the global Torelli theorem 
(\ref{_global_Torelli_Theorem_}),
we immediately obtain the following result.

\hfill

\theorem\label{_mono_ref_Theorem_}
Let $M=K3^{[n]}$ be a Hilbert scheme of points on K3,
${\cal M}_b$ its birational Teichm\"uller space,
and ${\cal M}_b(I)$ a connected component of ${\cal M}_b$.
Then ${\cal M}_b(I) \cong \Perspace/\Ref$,
where $\Perspace$ is the period domain defined
as in \eqref{_perspace_lines_Equation_},
and $\Ref=\Ref(H^2(M, \Z),q)$ the
corresponding reflection group,
acting on $\Perspace$ in a natural way.
\endproof

\hfill

The reflection group was computed in
 \cite{_Markman:constra_} (Lemma 4.2).
When $n-1$ is a prime power, this
computation is particularly effective.

\hfill

\definition\label{_spinorial_norm_Definition_}
Let $(V,g)$ be a real vector space
equipped with a non-de\-ge\-ne\-rate quadratic form
of signature $(m,n)$, and \[ S:=\{v\in V\ \ | \ \ g(v,v)>0\}.\]
It is easy to see that $S$ is homotopy equivalent
to a sphere $S^{m-1}$. Define the {\bf spinorial norm}
of $\eta\in O(V)$ as $\pm 1$, where the sign
is positive if $\eta$ acts as 1 on $H^{m-1}(S)$,
and negative if $\eta$ acts as -1.
Let $O^+(V)$ denote the set of
all isometries with spinorial norm 1.

\hfill

\remark 
It is easy to see that $\Ref\subset O^+(V)$,
where $\Ref$ is a reflection group.

\hfill

\proposition\label{_Refle_compu_Proposition_}
(\cite[Lemma 4.2]{_Markman:constra_}).
Let $M= K3^{[n]}$ be a Hilbert scheme of K3, 
and $\Ref=\Ref(H^2(M, \Z),q)$ the
corresponding reflection group. Then
$\Ref=O^+(H^2(M, \Z),q)$ if and only if
$n-1$ is a prime power.
\endproof

\hfill

\definition\label{_spin_orie_Definition_}
Let $V$ be a real vector space equipped
with a non-\-de\-ge\-ne\-rate quadratic form of signature 
$(m,n)$. A choice of {\bf spin orientation} on $V$
is a choice of a generator of the cohomology
group $H^{m-1}(S)$ (\ref{_spinorial_norm_Definition_}).
Clearly, $O^+(V)$ is a group of orthogonal maps
preserving the spin orientation.

\hfill

\remark
For a space $V$ with signature $(m,n)$, the group $O(V)$
has 4 connected components, which are given by a choice of orientation
and spin orientation. Alternatively, these 4 components are distinguished
by a choice of orientation on positive $m$-dimensional
planes and negative $n$-dimensional planes. 

\hfill

\remark
Donaldson (\cite{_Donaldson:polynomial_}) 
has shown that any 
diffeomorphism of a K3 surface $M$ preserves the spin orientation,
and the global Torelli theorem implies that every integer
isometry of $H^2(M)$ preserving the spin orientation
is induced by a diffeomorphism (\cite{_Borcea:Diffeomorphism_}).
This implies that the mapping class group $\Gamma_M$ is mapped
to $O^+(H^2(M, \Z))$ surjectively.

\hfill

\remark\label{_spin_orienta_hk_Remark_} 
Let $V= H^2(M, \R)$ be the second cohomology of
a hyperk\"ahler manifold, equipped with the Hodge structure 
and  the BBF form, and $V^{1,1}\subset V$ the space of
real (1,1)-classes. The set of vectors \[ R:=\{v\in V^{1,1}\ \
|\ \ q(v,v)>0\}\] is disconnected, and has two connected
components. Since the orthogonal complement
$(V^{1,1})^\bot$ is oriented, a spin orientation on $V$ is uniquely 
determined by a choice of one of two components of $R$. 
The K\"ahler cone of $M$ is contained in one of two
components of $R$. This gives 
a canonical spin orientation on $H^2(M, \R)$.

\hfill

\definition
Let $M$ be a hyperk\"ahler manifold. We say that {\bf the
Hodge-theoretic Torelli theorem holds 
for $M$}, if for any $I_1, I_2$ inducing 
 isomorphic Hodge structures on $H^2(M)$, the manifold
$(M,I_1)$ is bimeromorphically equivalent to 
$(M, I_2)$, provided that this isomorphism
of Hodge structures is also compatible with the
spin orientation and the 
Bogomolov-Beauville-Fujiki form, and $I_1, I_2$
lie in the same connected component of the
moduli space.

\hfill

\remark
This is the most standard version of 
global Torelli theorem.

\hfill

The following claim immediately follows from 
\ref{_mono_via_mappi_Theorem_}.

\hfill

\claim
 Let $M$ be a hyperk\"ahler manifold. 
Then the following statements are equivalent.
\begin{description}
\item[(i)] The Hodge-theoretic Torelli theorem holds 
for $M$.
\item[(ii)] The monodromy group $\Mon$ of $M$
is surjectively mapped to the group $O^+(H^2(M, \Z),q)$,
under the natural action of $\Mon$ on $H^2(M)$.
\end{description}
\endproof

\hfill

Comparing this with the Markman's computation
of the monodromy group (\ref{_Refle_compu_Proposition_}), 
we immediately obtain the following theorem.

\hfill

\theorem\label{_Torelli_prime_powe_Theorem_}
Let $M=K3^{[p^\alpha+1]}$.
Then the Hodge-theoretic Torelli theorem holds.
\endproof

\hfill

\remark
For other examples of hyperk\"ahler manifolds, the
Hodge-\-theo\-re\-tic global Torelli theorem is known to be false.
For some of generalized Kummer varieties this was proven
by Namikawa (\cite{_Namikawa:Torelli_}), and for
$M=K3^{[n]}$,
$n \neq p^\alpha+1$,
this observation is due to Markman (\cite{_Markman:constra_}).
For O'Grady's examples of hyperk\"ahler manifolds (\cite{_O_Grady_}),
it is unknown.


\section{Appendix: 
A criterion for a covering map (by Eyal Markman)}
\label{_Appe_Section_}


Another version of the proof of \ref{_exce_empty_covering_Proposition_}
was proposed by E. Markman; with his kind permission,
I include it here.

\hfill

\proposition (\ref{_exce_empty_covering_Proposition_})
\label{prop-criterion}
Let $\psi:X\rightarrow Y$ be a local homeomorphism
of  Hausdorff topological manifolds. 
Assume that every open subset $U\subset Y$, whose closure
$\overline{U}$ is homeomorphic
to a closed ball in $\R^n$, and such that $U$ is the interior of its closure,
satisfies the following property. For every connected component $C$
of $\psi^{-1}(\overline{U})$, the equality $\psi(C)=\overline{U}$ holds.
Then $\psi$ is a covering map.

\hfill

Verbitsky stated the above proposition 
in the category of differentiable manifolds and 
provided a proof of the proposition, 
involving Riemannian-geometric constructions on the domain $X$. 
We translate his proof to an elementary point set topology language. 
The natural translation of the statement and 
its proof to the category of differentiable  
manifolds is valid as well. In that case $\psi$ is a local diffeomorphism and 
it suffices for the assumption to hold 
for open subsets $U$, such that the boundary $\partial U$ is smooth,
and there exists a homeomorphism from $\overline{U}$ onto a closed ball in 
$\R^n$, which restricts to a diffeomorphism between the two interiors and 
between the two boundaries.
We will need the following well known fact (see \cite{_Browder_}, Lemma 1).

\hfill

\lemma
\label{lemma-liftingis-unique}
Let $f:X\rightarrow Y$ be a local homeomorphism of topological spaces,
$W$ a connected 
Hausdorff topological space, $h:W\rightarrow Y$ a continuous map, 
$x_0$ a point of $X$, and $w_0$ a point of $W$ satisfying $h(w_0)=f(x_0)$.
Then there exists at most one continuous map $\tilde{h}:W\rightarrow X$,
satisfying $\tilde{h}(w_0)=x_0$, and $f\circ\tilde{h}=h$.

\hfill

\noindent
{\bf Proof of \ref{prop-criterion}:}
The statement is local, so we may assume that $Y=\R^n$.
Let $x$ be a point of $X$ and set $y:=\psi(x)$.

\hfill

\definition
\label{def-x-star-shaped}
An open subset $U\subset \R^n$ is said to be $x$-{\em star-shaped}, 
if it satisfies the following conditions.
\begin{enumerate}
\item
$y$ belongs to $U$.
\item 
For every point $u\in U$, the line segment from 
$y$ to $u$ is contained in $U$. 
\item
There exists a continuous map
$\gamma:U\rightarrow X$, satisfying $\gamma(y)=x$,
and $\psi\circ\gamma:U\rightarrow \R^n$ is the inclusion.
\end{enumerate}

\claim
\label{claim-on-star-shaped-sets}
\begin{enumerate}
\item
\label{claim-item-finite-intersection}
Let $\{U_i\}_{i\in I}$ be a finite collection of
$x$-star-shaped open subsets of $\R^n$. Then their intersection
$\cap_{i\in I}U_i$ is $x$-star-shaped.
\item
\label{claim-item-union}
Let $\{U_i\}_{i\in I}$ be an arbitrary collection of
$x$-star-shaped open subsets of $\R^n$. Then their union
$U:=\cup_{i\in I}U_i$ is $x$-star-shaped.
\item
\label{claim-item-admissibility-criterion}
Let $U\subset \R^n$ be an $x$-star-shaped open subset,
$W\subset \R^n$ a connected open subset satisfying the following conditions.
a) $W\cap U$ is connected. b) For every point $t\in W\cup U$, the
line segment from $t$ to $y$ is contained in $W\cup U$. 
c) There exists a continuous map $\eta:W\rightarrow X$, such that
$\psi\circ\eta:W\rightarrow \R^n$ is the inclusion. d) There exists a point
$t\in W\cap U$, such that $\eta(t)=\gamma(t)$, where $\gamma:U\rightarrow X$
is the lift of the inclusion satisfying $\gamma(y)=x$. Then
$W\cup U$ is $x$-star-shaped.
\end{enumerate}

{\bf Proof:}
Part \ref{claim-item-finite-intersection} is clear.
Proof of part \ref{claim-item-union}: 
Let $\gamma_i:U_i\rightarrow X$
be the unique lift of the inclusion, satisfying
$\gamma_i(y)=x$. Define $\gamma:U\rightarrow X$ by
$\gamma(t)=\gamma_i(t)$, if $t$ belongs to $U_i$.
It sufficed to prove that $\gamma$ is well defined. 
If $t$ belongs to $U_i\cap U_j$, then $U_i\cap U_j$ is connected, being
$x$-star-shaped, and $\gamma_i(t)=\gamma_j(t)$, by
\ref{lemma-liftingis-unique}.

The proof of part \ref{claim-item-admissibility-criterion}
is similar to that of part \ref{claim-item-union}.
\endproof

\hfill

Given a positive real number $\epsilon$, set 
$B_\epsilon(y):=\{y'\in \R^n \ : \ d(y,y')<\epsilon\}$,
where $d(y',y)$ is the Eucleadian distance from $y'$ to $y$.
Let $\overline{B}_\epsilon(y)$ be the closure of $B_\epsilon(y)$.

\hfill

\claim
\label{key-claim}
Assume that $B_\epsilon(y)$ is $x$-star-shaped and let
$\gamma:B_\epsilon(y)\rightarrow X$ be the lift of the inclusion
satisfying $\gamma(y)=x$, as in  \ref{def-x-star-shaped}.
Then there exists an open connected subset $V\subset X$,
such that $V$ contains the closure $\overline{\gamma[B_\epsilon(y)]}$,  
$\psi:V\rightarrow \psi(V)$ is injective, and $\psi(V)$ is $x$-star-shaped. 

\hfill

{\bf Proof:}
Let $z$ be a point on the boundary $\partial\gamma[B_\epsilon(y)]$.
Then $\psi(z)$ belongs to the boundary of $B_\epsilon(y)$.
Now $\psi(z)$ has a basis of open neighborhoods $W$ with the property
that $U_z:=W\cup B_\epsilon(y)$ is $x$-star-shaped
(use \ref{claim-on-star-shaped-sets}
part \ref{claim-item-admissibility-criterion}). 
Let ${\cal U}_z$ denote the collection of all such $U_z$.
The collection 
$\{B_\epsilon(y)\}\cup 
\left[\cup_{z\in \partial \gamma[B_\epsilon(y)]}{\cal U}_z\right]$
is thus a collection of $x$-star-shaped open subsets.
Their union $U$ is $x$-star-shaped, by \ref{claim-on-star-shaped-sets},
so the inclusion $U\subset \R^n$ admits a lift $\gamma:U\rightarrow X$
satisfying $\gamma(y)=x$. Set $V:=\gamma[U]$.
Then $V$ is open, since $\gamma$ is a local-homeomorphism,
and $V$ contains the closure of $\gamma[B_\epsilon(y)]$,
by construction.
\endproof

\hfill

Let $D_x\subset \R^{>0}$  be the set of all $\epsilon\in \R^{>0}$,
such that there exists a continuous map 
$\gamma :\overline{B}_\epsilon(y)\rightarrow X$, satisfying $\gamma(y)=x$,
and such that $\psi\circ\gamma:\overline{B}_\epsilon(y)\rightarrow \R^n$ 
is the inclusion. 
Clearly, $D_x$ is a non-empty connected interval having $0$ as its left 
boundary point. We need to show that $D_x=\R^{>0}$.
It suffices to show that $D_x$ is both open and closed.

\hfill

\claim
\label{claim-open}
$D_x$ is open.

\hfill

{\bf Proof:}
Let $\epsilon$ be a point of $D_x$. The image 
$\gamma[\overline{B}_\epsilon(y)]$ is compact and $X$ is Hausdorff.
Hence, $\gamma[\overline{B}_\epsilon(y)]$ is closed and is thus equal to the 
closure of $\gamma[B_\epsilon(y)]$.
Then $\psi\left(\overline{\gamma[B_\epsilon(y)]}\right)=\overline{B}_\epsilon(y)$.
Hence, there exists an open $x$-star-shaped subset $U\subset \R^n$,
containing $\overline{B}_\epsilon(y)$, by \ref{key-claim}.
Compactness of $\overline{B}_\epsilon(y)$ implies that $U$
contains $\overline{B}_{\epsilon_1}(y)$, for some $\epsilon_1>\epsilon$. 
Now $\epsilon_1$ belongs to $D_x$, since $U$ is $x$-star-shaped.
Hence, $D_x$ is open.
\endproof

\hfill

Set $s:=\sup(D_x)$. If $s$ is infinite, we are done. Assume that $s$ is finite.
$B_s(y)$ is $x$-star-shaped, by \ref{claim-on-star-shaped-sets}.
Let $\gamma:B_s(y)\rightarrow X$ be the lift of the inclusion satisfying
$\gamma(y)=x$.

\hfill

\claim
\label{claim-closure}
The closure $C:=\overline{\gamma[B_s(y)]}$ is a connected 
component of the preimage $\psi^{-1}[\overline{B}_s(y)]$. Furthermore,
$\psi:C\rightarrow \overline{B}_s(y)$
is injective.

\hfill

{\bf Proof:}
There exists an open subset $V$ of $X$, containing $C$, 
such that $\psi:V\rightarrow \psi(V)$ is 
a homeomorphism, by \ref{key-claim}. Hence, 
$V\cap \psi^{-1}[\overline{B}_s(y)]=C$, and $C$ is both open and closed in 
$\psi^{-1}[\overline{B}_s(y)]$.
\endproof

\hfill

Up to now we used only the assumption that $\psi$ is a local homeomorphism.
We now use the assumption that $\psi:C\rightarrow \overline{B}_s(y)$
is surjective, for every connected component of $\psi^{-1}[\overline{B}_s(y)]$, 
and in particular for $C:=\overline{\gamma[B_s(y)]}$. 
We conclude that $s$ belongs to $D_x$.
A contradiction, since $D_x$ is open.
This completes the proof of \ref{prop-criterion}.
\endproof


\section{Appendix: Subtwistor metric on the period space}
\label{_Appe_subtwistor_}


\subsection{GHK lines and subtwistor metrics}

Let $\Perspace$
be a period space of a hyperk\"ahler manifold $M$, 
identified with a Grassmannian $SO(b_2-3,3)/SO(2)\times SO(b_2-1,3)$
of oriented, positive 2-planes in $H^2(M,\R)$.
We shall consider $\Perspace$ as a complex manifold,
with the complex structure obtained as in \eqref{_perspace_lines_Equation_}.
Fix an auxiliary Euclidean 
metric $g$ on $H^2(M)$. Given a positive 3-dimensional
plane $W\subset H^2(M)$, denote by $S_W\subset \Perspace$
the set of all 2-dimensional oriented planes contained in $W$.
Clearly, $S_W$ is a complex curve in $\Perspace$.
The metric $g$ induces a Fubini-Study metric on $S_W$.

Consider a sequence $S_1, ..., S_n$ of intersecting hyperk\"ahler 
lines connecting $x\in \Perspace$ to $y\in \Perspace$,
with $s_i \in S_i \cap S_{i+1}$, $i = 1, ..., n-1$
the intersection points, and $s_0:=x$, $s_{n+1}:=y$.
Denote by $l_{S_1, ..., S_n}(x,y)$ the sum $\sum_{i=0}^n d_g(s_i, s_{i+1})$,
where the distance $d_g(s_i, s_{i+1})$ is computed on the
hyperk\"ahler line $S_{i+1}$ using the metric induced by $g$
as above.  Let 
\[ d_{tw}(x,y):= \inf_{S_1, ..., S_n}l_{S_1, ..., S_n}(x,y)
\]
where the infimum is taken over all appropriate 
sequences of GHK lines, connecting $x$ to $y$.

The following theorem is stated for periods
of hyperk\"ahler manifolds, but in fact it could
be stated abstractly for $\Perspace=SO(m-3,3)/SO(2)\times SO(m-1,3)$,
for any $m>0$. No results of geometry or topology of
hyperk\"ahler manifolds are used in its proof.

\hfill

\theorem\label{_d_tw_metric_Theorem_}
Let $\Perspace$ be a period space of a hyperk\"ahler manifold,
and $d_{tw}:\; \Perspace\times\Perspace \arrow \R^{\geq 0}\cup \infty$
the function defined above. Then
\begin{description}
\item[(i)] $d_{tw}$ is a metric on $\Perspace$
\item[(ii)] The metric $d_{tw}$ induces the
usual topology on $\Perspace$.
\end{description}
The rest of this section is
taken by the proof of \ref{_d_tw_metric_Theorem_}.

\hfill

\definition
The metric $d_{tw}$ is called {\bf the subtwistor metric}
on the period space, and a piecewise geodesic connecting $x$ 
to $y$ and going over $S_i$ is called {\bf a subtwistor path}.

\hfill

\remark
It is in many ways similar to the sub-Riemannian
metrics known in metric geometry (see e.g. \cite{_BBI_}).

\hfill

The triangle inequality for $d_{tw}$ is clear from its definition.
To prove that $d_{tw}$ is a metric, we need only to show that
$d_{tw}< \infty$ and $d_{tw}(x,y)>0$ for $x\neq y$.
The inequality $d_{tw}< \infty$ follows from \ref{_4_GHK_lines_Proposition_},
which claims that all points of $\Perspace$ are connected
by a finite sequence of GHK lines.
The latter condition is clear, because
$d_{tw} \geq d_g$, where $d_g$ is a geodesic
distance function on $\Perspace$ associated with the
Riemannian metric $g$.

\subsection{Subtwistor metric on a Lie group}

Let $(\Perspace, d_g)$ be the space $\Perspace$
equipped with a Riemannian metric associated with 
a scalar product $g$ on $V:=H^2(M,\R)$, and $d_{tw}$
its subtwistor metric.
To finish the proof of \ref{_d_tw_metric_Theorem_},
we have to show that the identity map
$(\Perspace, d_{tw}) \stackrel \tau \arrow (\Perspace, d_g)$
is a homeomorphism. 

Notice that $\tau$ is continuous, because $d_{tw}\geq d_g$.
Brouwer's
invariance of domain theorem implies that
it suffices to show that $(\Perspace, d_{tw})$
is homeomorphic to a manifold.

\hfill

\claim\label{_bije_cont_homeo_Claim_}
(Brouwer's invariance of domain theorem)
Let $X\stackrel f\arrow Y$ be a continuous, bijective
map of Hausdorff manifolds. Then $f$ is a homeomorphism.

\hfill

{\bf Proof:} This is a corollary of L. E. J. Brouwer's Theorem
on invariance of domain, proven in Beweis der Invarianz des n-dimensionale
Gebiets, Math. Annalen 71 (1911), pages 305-315. See also Terence Tao's blog:
{\tiny \url{http://terrytao.wordpress.com/2011/06/13/brouwers-fixed-point-and-invariance-of-domain-theorems-and-hilberts-fifth-problem/}}\footnote{I am grateful to the referee for this observation and the reference.}
 \endproof

\hfill

We prove that $(\Perspace, d_{tw})$ is a manifold by an application
of the Gleason-Palais theorem on transformation groups,
obtained in 1950-ies
as a byproduct of the solution of Hilbert's 5th problem
by Gleason, Montgomery and Zippin.
This would require us to switch from $\Perspace$ to a 
Lie group $G\subset SO(H^2(M, \R), q)$ transitively acting on $\Perspace$. 
We introduce a metric $d_{tw}$ on $G$, in such a way
that $(\Perspace, d_{tw})$ is obtained as a quotient
of $(G, d_{tw})$, and prove that $(G, d_{tw})$ is a manifold,
using the Gleason-Palais theorem.

\hfill

Let $V:=H^2(M, \R)$, considered as a vector space
with the scalar product $q$ of signature $(3, n-3)$, 
and let $G$ be the connected component of identity of $SO(V)$.

\hfill

\definition
An {\bf elementary transform} is an element $h\in G$ 
fixing a codimension 2 subspace $V_1\subset V$ of
signature $(1,n-3)$. {\bf An elementary decomposition}
of $h\in G$ is a decomposition $h=h_1h_2...h_n$,
where $h_i$ are elementary transforms.

\hfill

\remark
Any element of $G$ admits an elementary decomposition,
obviously non-unique. This is proven by the same argument
as used to show that $SO(n)$ is generated by rotations
fixing a codimension 2 subspace.

\hfill

For any elementary transform $h\in G$, 
choose orthogonal coordinates in which $h$
is a turn, with an angle $0\leq \alpha \leq \pi$,
and denote by $\|h\|$ the number $|\alpha|$.
Consider the Lie algebra element $\log h$,
corresponding to $h$, and let ${\goth h}\in T\Perspace$
be the corresponding tangent vector field. For any
$x\in \Perspace$, one has $T_x\Perspace = \Hom(x, x^\bot)$,
where $x^\bot$ is an orthogonal complement to $x$
which is considered as a codimension 2 subspace in $V$.
Fix an auxiliary positive definite Euclidean metric $g$ on $V$.
Clearly, at $x\in \Perspace$, one has
\begin{equation}\label{_length_of_ta_ve_Equation_}
|{\goth h}|_g\restrict x = \left|\pi_{x^\bot}\left(\log h \restrict x\right)\right|_g
\end{equation}
where $\pi_{x^\bot}$ is an orthogonal (with respect to
$g$)  projection to $x^\bot$. Since $x^\bot$ is
2-dimensional, the quantity \eqref{_length_of_ta_ve_Equation_}
is bounded by $|a_1|+|a_2|$, where $a_i$ are eigenvalues 
of $\log h$. This gives
\begin{equation}\label{_length_of_ta_ve_via_norm_Equation_}
|{\goth h}|_g\restrict x \leq 2 \|h\|.
\end{equation}

Consider the path
 $\gamma:\; [0,1] \arrow \Perspace$, $\gamma(t)=
e^{t\log h}(x)$, connecting $x$ to $y:=h(x)$.
Then $d_g(x, y) \leq \int_0^1 |{\goth h}(e^{t\log h}(x))| dt$.
Therefore, \eqref{_length_of_ta_ve_via_norm_Equation_}
gives $d_g(x, h(x)) \leq 2 \|h\|$.

\hfill

\definition
Define the {\bf subtwistor norm} on $G$
as $\|h\|_{tw}:= \inf(\|h_1\|+\|h_2\|+...+\|h_n\|)$,
where the infinum is taken over all elementary
decompositions $h=h_1 h_2 ... h_n$.

\hfill

\remark 
It is easy to check that this norm satisfies the usual axioms,
that is, defines a right-invariant metric on the Lie group,
using the formula $d_{tw}(x,y):=\|xy^{-1}\|_{tw}$.

\hfill

\claim\label{_subtw_norm_group_perspace_Claim_}
Let $x,y \in \Perspace$. Then
$\mu\leq d_{tw}(x,y)\leq 2\mu$, where 
\[ \mu:= \inf_{h\in G, h(x)=y}\|h\|_{tw},
\]
and infimum is taken over all $h\in G$ such that $h(x)=y$.

\hfill

{\bf Proof:} Let $S_1, ..., S_n$ be a sequence of 
GHK lines connecting $x$ to $y$, and $x_0, ..., x_n$ 
the intersection points with $x_0 =x, x_i \in S_{i}\cap S_{i+1}$
and $x_n =y$. By definition, $d_{tw}(x,y)$ is an infimum of
$\sum_i d_g(x_i, x_{i+1})$ for all such sequences.
Let $W_i\subset V$ be a 3-dimensional positive 
subspace corresponding to $S_i$, and $h_i$ an elementary
transform acting trivially on $W_i^\bot$ and
mapping $x_{i-1}$ to $x_i$ (such $h_i$ exists, because
$SO(3)$ acts transitively on 2-planes).\footnote{The corresponding
rotation can be chosen in such a way that its angle is
equal to the distance between $S_{i-1}$ and $S_i$ in the twistor
line; then $\|h_i\|\geq d_g(x_{i-1},x_i)$.}

 Then $h:= h_1... h_n$
maps $x$ to $y$. The number $d_g(x_{i-1},x_i)$
is equal to $d_{tw}(x_{i-1},x_i)$, because these two points
lie on the same twistor line. Moreover, this number
is equal to a length of the smallest circle
segment in $S_{i}$ connecting $x_{i-1}$ to $x_i$,
which is equal to $\|h_i\|$. We obtained
\[
d_{tw}(x,y) = \inf_{\{S_i\}} \sum_i d_g(x_{i-1},x_i)=
\sum_i \|h_i\| \geq \|h\|_{tw}\geq \mu.
\]
Using \eqref{_length_of_ta_ve_via_norm_Equation_},
we also obtain 
\[
d_{tw}(x,y) \leq \sum_i d_g(x_{i-1},x_i)\leq  \sum_i 2 \|h_i\| = 2\mu.
\]
\endproof

\hfill

From \ref{_subtw_norm_group_perspace_Claim_},
the following observation is apparent.

\hfill

\corollary\label{_dista_between_orbits_subtwi_Corollary_}
Let $G\subset SO(V)$ be a connected component, $\Perspace =Gr_{++}(V)$
the period space, $x\in \Perspace$ a point, 
and $G\stackrel {\tau_x}\arrow \Perspace$ the map mapping $g$ to $g(x)$. 
Let $y\in \Perspace$, and let $A$ be the distance 
between $\tau_x^{-1}(x)$ and $\tau_x^{-1}(y)$, computed
with respect to the subtwistor norm on $G$. Then 
$A \leq d_{tw}(x,y)\leq 2A$.

\hfill

{\bf Proof:} Clearly, $A$ is an infimum of
$\|vu^{-1}\|_{tw}$ for all $u,v\in G$ satisfying
$u(x)=x, v(x)=y$. In terms of
\ref{_subtw_norm_group_perspace_Claim_},
this number is equal to $\mu$, and then
the inequality $A \leq d_{tw}(x,y)\leq 2A$
follows. \endproof

\hfill

Given a norm (or a right-invariant metric) on a group,
one can construct a metric on its right quotients, as follows.
One defines the
metric on the space of right classes $G/G_x$, 
using the Hausdorff distance between 
the right classes $yG_x$ and $zG_x$
computed with respect to the given metric on 
$G$.

\hfill

\corollary\label{_quotie_subtwi_Corollary_}
Let $G\subset SO(V)$ be a subgroup defined above,
$\|\cdot \|$ the subtwistor norm, $\Perspace =Gr_{++}(V)$
the corresponding period space, $x,y,z\in \Perspace$, 
and $G\stackrel {\tau_x}\arrow \Perspace$ the map mapping
$g$ to $g(x)$. Denote by $G_x$ the stabilizer of $x$ in
$G$, $G_x:= \tau_x^{-1}(x)$. We equip $G/G_x$ with
a metric, using the subtwistor metric on $G$ as above. 
 Then the natural map $\tau_x:\; G \arrow \Perspace$
induces a homeomorphism from $G/G_x$ to $\Perspace$.

\hfill

{\bf Proof:}
Let $A_{x,y}$ be the Hausdorff distance between 
$\tau_x^{-1}(y)$ and $\tau_x^{-1}(z)$, computed in 
the subtwistor metric. To prove
\ref{_quotie_subtwi_Corollary_}, it would suffice to show
that for some number $\mu(z)>1$ depending on $z$, one has
\begin{equation}\label{_mu_of_z_Equation_}
\mu(z)^{-1}A_{x,y}\leq d_{tw}(y,z)\leq \mu(z) A_{x,y}
\end{equation}
 When $z=x$,
this inequality is implied directly
by \ref{_dista_between_orbits_subtwi_Corollary_}:
 \[ A_{x,y} \leq d_{tw}(y, x)\leq 2A_{x,y}.
\]
For general $z$, we use the action of $G$ on
itself and $\Perspace$.

The group $G:=SO(H^2(M, \R), q)$
acts on $\Perspace$ transitively, and each $\gamma\in G$
induces a bi-Lipschitz map on $\Perspace$, distorting the metric 
$d_g$ in a way which is bounded by $C_\gamma$, 
where $C_\gamma$ is a constant depending on
the largest eigenvalue of $\gamma(g)g^{-1}$.
Indeed, $d_g$ is the metric on the Grassmannian
of 2-planes in $H^2(M,\R)$ associated with the metric $g$
on $H^2(M,\R)$, and $\gamma$ distorts this metric 
in a controlled way, with the Lipschitz constant
bounded by the eigenvalues of $\gamma(g^{-1}) g$.

The same argument shows that $\gamma$
distorts $(\Perspace, d_{tw})$ by the same
number: 
\[ 
C_\gamma^{-1}d_{tw}(x,y)\leq   d_{tw}(\gamma(x), \gamma(y))\leq C_\gamma
d_{tw}(x,y).
\]
Now, let $\gamma\in G$ be an element mapping $z$ to $x$.
Then the right action of $\gamma$ maps $\tau_x$ to $\tau_z$
and the pair $(y,z)$ to $(\gamma(y), x)$, giving
\begin{equation}\label{_A_dista_gamma_Equation_}
A_{x,y} \leq d_{tw}(\gamma(y), x)\leq 2A_{x,y} 
\end{equation}
by 
the above argument. Since the action
of $\gamma$ on $\Perspace$ is bi-Lipschitz,
\begin{equation}\label{_bi-Lipschi_gamma_Equation_}
C_\gamma^{-1}d_{tw}(y,z) \leq d_{tw}(\gamma(y), x)\leq C_\gamma d_{tw}(y,z).
\end{equation}
Comparing this inequality with \eqref{_A_dista_gamma_Equation_}
we obtain \eqref{_mu_of_z_Equation_}:
\[
\frac{A_{x,y}}{2C_\gamma} \leq \frac{A_{x,y}}{C_\gamma}
\stackrel{\text{\scriptsize\eqref{_A_dista_gamma_Equation_}}} \leq
\frac{d_{tw}(\gamma(y), x)}{C_\gamma} 
\stackrel{\text{\scriptsize\eqref{_bi-Lipschi_gamma_Equation_}}}\leq d_{tw}(y,z)
\stackrel{\text{\scriptsize\eqref{_bi-Lipschi_gamma_Equation_}}}\leq 
C_\gamma d_{tw}(\gamma(y), x) 
\stackrel{\text{\scriptsize\eqref{_A_dista_gamma_Equation_}}} \leq
2C_\gamma A_{x,y}.
\]
\endproof

\hfill

\corollary\label{_quotie_obtain_subtwi_Corollary_}
Let $G\subset SO(H^2(M, \R), q)$ be the connected component of unit
in $SO(H^2(M, \R), q)$ 
acting on the period space $\Perspace$, and $d_{tw}$
the metric defined from the subtwistor norm. 
Then $(\Perspace, d_{tw})$ is homeomorphic to
the quotient $G/G_x$, equipped with 
the quotient topology induced from $d_{tw}$.

{\bf Proof:} Follows from 
\ref{_quotie_subtwi_Corollary_}. \endproof

\hfill

To prove that $(\Perspace, d_{tw})= (G,d_{tw})/G_x$ is a manifold,
it would suffice to show that $G$ with the metric
induced from a subtwistor norm is a manifold.
Indeed, in this case, the subtwistor norm
induces the standard topology on $G$
by \ref{_bije_cont_homeo_Claim_}. Then
$(\Perspace, d_{tw})$ is a manifold
by \ref{_quotie_obtain_subtwi_Corollary_}.

\hfill

\theorem\label{_grp_w_subtwi_mfld_Theorem_}
Let $G\subset SO(H^2(M, \R), q)$ be the Lie group
defined as above, and $d_{tw}$
the metric defined from the subtwistor norm. 
Then $(G,d_{tw})$ is a topological manifold.

\hfill

We prove \ref{_grp_w_subtwi_mfld_Theorem_}
in Subsection \ref{_Glea_Pala_Subsection_}.

\subsection{Gleason-Palais theorem and its applications}
\label{_Glea_Pala_Subsection_}

\ref{_grp_w_subtwi_mfld_Theorem_}
follows directly from a theorem of
Gleason and Palais about 
transformation groups (\cite{_Gleason_Palais_};
for a more recent treatment and reference, see
\cite{_Banakh_Zdomsky_}).

\hfill

\definition
Let $M$ be a topological space. We say that $M$ {\bf has
  Lebesgue covering dimension $\leq n$} if every
open covering of $M$ has a refinement $\{U_i\}$
such that each point of $M$ belongs to at most 
$n+1$ element of $\{U_i\}$. A {\bf Lebesgue covering
dimension} of $M$ (denoted by $\dim M$) is an infimum
of all such $n$. 

\hfill

The following two well-known claims are easy to prove.

\hfill

\claim\label{_dim_mfld_Claim_}
If $M$ is an $n$-manifold, $\dim M=n$. \endproof

\hfill

\claim\label{_dim_subset_Claim_}
If $X\subset M$ is a subset of a
topological space, with induced topology, 
one has $\dim X \leq \dim M$. \endproof

\hfill

The following theorem is a deep and important result,
obtained in 1950-ies in the course of studying transformation groups
along with the solution of Hilbert's 5-th problem.

\hfill

\theorem\label{_Gleason_Palais_Theorem_}
(Gleason-Palais)
Let $G$ be a topological group, which is
locally path connected, and has $\dim K < \infty$ for
each compact, metrizable subset $K\subset G$. Then  $G$ is
homeomorphic to a Lie group. 

\hfill

{\bf Proof:} \cite[Theorem 7.2]{_Gleason_Palais_}.
\endproof

\hfill

Now we can finish the proof of 
\ref{_grp_w_subtwi_mfld_Theorem_}, and 
\ref{_d_tw_metric_Theorem_}. The group $(G, d_{tw})$
is by construction locally path connected. Moreover, one has 
$d_{tw}\geq d_g$, where $d_g$ is a metric obtained 
from a positive definite metric on $V$; this is 
proven in the same way as the inequality $d_{tw}\geq d_g$
for $\Perspace$.  Therefore, the
map $(G, d_{tw})\arrow (G, d_g)$ is continuous.
This implies that any compact $K\subset (G, d_{tw})$
is homeomorphic to its image in $(G,d_g)$, hence it is 
finitely-dimensional. Applying Gleason-Palais, we
obtain that $(G, d_{tw})$ is a manifold.
\endproof

\hfill

{\bf Acknowledgements:}
This paper owes much to Eyal Markman, whose efforts
to clean up my arguments were invaluable.
Many thanks to E. Loojienga and F. Bogomolov
for their remarks and interest. An early version of this 
paper was used as a source of a mini-series of 
lectures at a conference ``Holomorphically symplectic
varieties and moduli spaces'', in Lille,
June 2-6, 2009. I am grateful to the organizers
for this opportunity and to the participants
for their insight and many useful comments.
This paper was written at Institute
of Physics and Mathematics of the Universe 
(University of Tokyo), and I would like to thank
my colleagues and the staff at IPMU 
for their warmth and hospitality. Many thanks
to Eyal Markman for many excellent suggestions, remarks and
corrections to several earlier versions of the manuscript, 
to Dan Huybrechts for his interest and comments,
and to Brendan Hassett and Yuri Tschinkel for finding gaps
in the argument.

{

\noindent {\sc Misha Verbitsky\\
{\sc Laboratory of Algebraic Geometry,\\
National Research University HSE,\\
Faculty of Mathematics, 7 Vavilova Str. Moscow, Russia,}\\
\tt  verbit@mccme.ru}, also: \\
{\sc Kavli IPMU (WPI), the University of Tokyo}

}

\end{document}